%% file: manuscript-calcolo-R2.tex
\IfFileExists{./prepreamble-calcolo.sty}{\RequirePackage[packages,theorems]{prepreamble-calcolo}}{}

\documentclass[%
sn-mathphys-num,
]{sn-jnl-scoop}

\IfFileExists{./preamble-calcolo.sty}{\RequirePackage[packages,theorems]{preamble-calcolo}}{}

\usepackage{manuscript}
\usepackage{KMS-local}
\usepackage{amsmath}
\usepackage{scoop-semantic}
\usepackage{scoop-syntax}
\usepackage{scoop-packages}
\usepackage[utf8]{inputenc}
\usepackage{KMS-standard-packages}

\IfFileExists{./postpreamble-calcolo.sty}{\RequirePackage[packages,theorems]{postpreamble-calcolo}}{}

\usepackage{manyfoot}
\usepackage{xcolor}
\usepackage{amsmath}

\makeatletter
\@ifpackageloaded{hyperref}{%
	\hypersetup{
		pdftitle = {A block Recycled GMRES method with investigations into aspects of solver performance},
		pdfauthor = {Michael L Parks, Kirk M Soodhalter, Daniel B Szyld},
		pdfkeywords = {Krylov subspace methods, deflation, subspace recycling, block Krylov methods, high-performance computing},
		pdfcreator = {Created using the Scoop Template Engine version 1.6.0.post0+git.eeaa8072.dirty.}
	}
}{
	\pdfinfo{
		/Title (A block Recycled GMRES method with investigations into aspects of solver performance)
		/Author (Michael L Parks, Kirk M Soodhalter, Daniel B Szyld)
		/Subject ()
		/Keywords (Krylov subspace methods, deflation, subspace recycling, block Krylov methods, high-performance computing)
		/Creator (Created using the Scoop Template Engine version 1.6.0.post0+git.eeaa8072.dirty.)
	}
}
\makeatother

\begin{document}

\title[Solver performance for block recycled GMRES]{A block Recycled GMRES method with investigations into aspects of solver performance}
\subtitle{}

\author[1]{Michael L Parks}
\email{\detokenize{parksml@ornl.gov}}
\author*[2]{Kirk M Soodhalter}
\email{\detokenize{ksoodha@maths.tcd.ie}}
\author[3]{Daniel B Szyld}
\email{\detokenize{szyld@temple.edu}}
\affil[1]{Computer Science and Mathematics Division, Oak Ridge National Laboratory, Oak Ridge, TN, 37830, USA}
\affil[2]{School of Mathematics, Trinity College Dublin, College Green, Dublin 2, Ireland}
\affil[3]{Department of Mathematics, Temple University, Philadelphia, PA, 19147, USA}

\abstract{%
\input{abstract-R2.tex}
}




\pacs[Funding]{The second author was supported in part by Research Ireland grand 22/EPSRC/3857.  The third author was supported in part by the U.S.\ National Science Foundation under grant DMS-1418882.}

\maketitle

\input{content-R2.tex}

\begin{appendices}

\end{appendices}

\IfFileExists{./manuscript-calcolo.bib}{\bibliography{manuscript-calcolo.bib}}{\bibliography{manuscript.bib}}

\end{document}

%% file: abstract-R2.tex
We propose a block Krylov subspace version of the \gcrodr method proposed in
[Parks et al.; SISC 2005], which is an iterative method allowing for the
efficient minimization of the residual over an augmented Krylov subspace.
We offer a clean derivation of our proposed method and discuss methods of
selecting recycling subspaces at restart as well as implementation decisions in
the context of high-performance computing.  Two types of numerical experiments are
presented: those demonstrating convergence properties and those demonstrating the
data movement and cache efficiencies of the dominant operations of the method,
measured using processor monitoring code from Intel.

%% file: content-R2.tex
\section{Introduction}
\label{section:introduction}
We explore the efficient solution of a sequence of linear systems
\begin{align}
	\bA^{[i]}
	\paren[auto]{(}{)}
	{
		\bX_0^{[i]}
		+
		\bT^{[i]}
	}
	=
	\bB^{[i]},
	\label{eqn.Axb}
\end{align}
where the coefficient matrices $\bA^{[i]}\in\C^{n\times n}$ are assumed to be non-Hermitian and the right-hand sides $\bB^{[i]}\in
\C^{n \times p}$ may or may not also change with $i$. For many applications, the matrices $\bA^{[i]}$ are large and sparse. The
block $\bX_0^{[i]}\in\C^{n\times p}$ is the initial approximations to the solution of system $i$ and $\bT^{[i]}$ is the
corresponding initial error.

At the core of many problems in the computational sciences is the need to solve large, sparse linear systems.  It is often the
case that one must solve a sequence of systems, and these systems are somehow related. Examples include: uncertainty quantification
\cite{JC-recyc-uncer,JC-dd-recyc-uncer.2007}, Newton-like iterations from, e.g., a density functional theory computations
\cite{Heroux2007}, topology optimization \cite{WSG.2007}, the modeling of crack propagation in materials~\cite{Parks.deSturler.GCRODR.2005}, and tomography \cite{Kilmer.deSturler.tomography.2006}. This challenge often includes solving
for more than one right-hand side for each $i$. For each coefficient matrix, one could simply treat each right-hand side in
sequence. However, it is often more efficient to take advantage of the underlying structure of all problems combined. Block Krylov
subspace iterative methods \cite{OLeary1980,Vital1990} were proposed to solve systems with many right-hand sides, and such methods
were proposed to accelerate convergence even when there is is only one right-hand side; see, e.g.,
\cite{CK.2010,OLeary1980,O-ParCG.1987,S.2014,Schaerer.Szyld.Torres.23}. Subspace recycling techniques \cite{Parks.deSturler.GCRODR.2005} were proposed to
take advantage of relationships between sequences of coefficient matrices.  It is therefore natural to combine these two
strategies to take advantage of their individual benefits as well as any synergistic interactions that arise.

This paper is the journal version of a technical report published in 2016 \cite{ParksSoodhalterSzyld:2016:1},
which was derived from work in the PhD thesis \cite{S.2012-thesis}.
We have polished and slightly updated the content, and added some numerical experiments.
We describe the original development of the combination of these
two techniques (block Krylov subspace methods and recycling), which led to implementations in \trilinos \cite{PS.GCRODR-Code.2011}
in 2011 and an accompanying \matlab implementation \cite{GCRODR-matlab-code.2016}.  We justify the utility of the resulting method
not just with steeper convergence curves but also with performance experiments. Specifically, we measure the
run-times and cache efficiency of key computational kernels to demonstrate the utility of such block methods.  
We mention that the recent paper \cite{Thomas.Baker.Gaudreault.SISC.25} complements our work. While we concentrate
on the implementation of a block version of GCRO-DR, on
reducing data movement, and on using cache more efficiently, the authors of 
\cite{Thomas.Baker.Gaudreault.SISC.25}
concentrate on the block orthogonalization step.

The rest of this paper is organized as follows. In \Cref{section.background},
we describe the problem being solved and discuss what has been investigated in
the literature. In \Cref{section.prelim}, we briefly review Krylov subspace
methods and their generalization to the block setting as well as a framework
for understanding subspace augmentation methods, including those employing a
recycling strategy.  In \Cref{section.deriv-impl}, we extend the recycled
\gmres method to the block Krylov subspace setting.  We further describe
implementation decisions meant to improve the data movement efficiency of the
method.  In \Cref{section.conv}, we discuss convergence properties and theory
of this method.  In \Cref{section.num}, we present experiments. These are
divided into two types.  The first includes simple convergence experiments,
demonstrating the competitiveness of these methods for large-scale problems. We
then present measurements of how data is moved to and used on the processor for
the dominant operations of the method.

\section{Background}\label{section.background}
Krylov subspace iterative methods are a standard tool for solving sparse systems such as those arising in \eqref{eqn.Axb}.  In
this work, we consider solving \eqref{eqn.Axb} using block Krylov subspace techniques (\Cref{subsection.block-krylov})
in the case that the number of right-hand sides $p>1$, but we also explore the
utility of using these techniques for the case that $p=1$, in order
to accelerate convergence. 
The latter application of the block techniques was first
suggested in \cite{OLeary1980} and elaborated upon in \cite{O-ParCG.1987} and also used in, e.g., \cite{CK.2010,S.2014};  see,
e.g.,~\cite{Gutknecht2007}, for a nice introduction to the topic.

We describe the original development of the combination of
block \gmres with the recycling-based augmentation scheme
\gcrodr, and the implementation decisions that were taken in the design of the codes available in 
\cite{PS.GCRODR-Code.2011,GCRODR-matlab-code.2016}.
In the case $p=1$, \gcrodr was introduced \cite{Parks.deSturler.GCRODR.2005} for the treatment of sequences of
\quotes{slowly-changing} linear systems. The expression slowly-changing is intentionally imprecise; it can mean that each matrix is a
small Frobenius norm perturbation of its predecessor, e.g., arising from during a Newton iteration or from the modeling of crack
propagation \cite{Littlewood.2016}. It can also refer to a sequence of systems whose spectral structure has some relationship
(though their norm distance from one another is nontrivial), as may be the case when evaluating multiple parameter realizations of
some underlying PDE model, e.g., in stochastic PDE applications \cite{JC-recyc-uncer}. This method allows one to retain important
approximate invariant subspace information generated during the solution of the $i$th linear system, and leverage that information
to accelerate convergence of the iteration to solve the subsequent $(i+1)st$ system.  The \gcrodr method is one of many subspace
augmentation-based recycling approaches.  A general framework can be found in \cite{SoodhalterDeSturlerKilmer:2020:1}, which
greatly simplifies the presentation of such approaches and enables more straightforward development of new techniques.

 As high-performance computing architectures continue to evolve, the cost of floating point computations has decreased
 dramatically when compared to the cost of data movement; this effects both algorithm performance and power consumption costs.
 See, e.g., \cite{Shalf.Dosanjh.Morrison.11} for further discussion and \cite{Boman.Higgins.Szyld.23} for
 an analysis of how the block size influences the performance. 
 Metrics such as the amount of data moved and the efficiency
 of cache reuse have become more important measures of algorithm performance than simply counting floating point operations; see,
 e.g., \cite{H.2010,MHDY.2009}.  \emph{Arithmetic intensity} (i.e., the amount of computation done per unit of memory accessed) is
 an effective quantification of how well a particular algorithm can perform in the HPC context.

 For example, in the dense linear algebra setting, it was shown that level-3 \blas matrix-matrix operations (such as multiplying a
 dense matrix times a block of vectors) demonstrate superior performance over level-2 \blas matrix-vector operations when measured in
 terms of arithmetic intensity. In the sparse linear algebra setting, the dominant operation of most Krylov subspace methods
 for the case $p = 1$ is a sparse matrix-vector multiplication, and in block Krylov subspace methods this is replaced with a
 sparse matrix-matrix multiplication (which in this context means a sparse matrix multiplied times a block of vectors, which is
 often dense). It was shown that this sparse block operation also demonstrates superior performance to its non-block
 counterparts in data-related metrics \cite{BDJ.2006,H.2010}. It is known from theory and observation that block Krylov methods
 converge in fewer iterations than their single-vector counterparts (and at a minimum can do no worse).  The superiority of
 performance when applying a sparse matrix to a block of $p$ vectors when compared to $p$ single-vector matrix-vector products (in
 terms of time and data-movement measurements) is clear, especially on modern and emerging architectures; 
 \cf~\cite{Boman.Higgins.Szyld.23} and
 \Cref{section.small-experiments}.

Given this observation, it is reasonable to consider a block version of \gmres with recycling to leverage both the algorithmic and
hardware advantages, and to explore their application to block systems and systems with a single right-hand side. The extension of
\gcrodr and other such methods to the block setting is a natural one to make.  This paper describes the development of the block \gcrodr
high-performance implementation \cite{PS.GCRODR-Code.2011} in the \belos package of the \trilinos project \cite{2011a}. Based upon
this code, the authors of \cite{GiraudJingXiang:2022:1} have extended this method to the flexible preconditioning setting and
treat the issue of inexact block Krylov subspace breakdown thoroughly.  Other such methods have also been extended to the block
setting; see, e.g., \cite{Darnell2008,MZL.2014}.  Based on the original implementation in \trilinos \cite{PS.GCRODR-Code.2011},
other authors have refined the block \gcrodr algorithm to make the HPC implementation more effective; see, e.g.,~\cite{AudibertGirardonHaddarJolivet:2023:1,JolivetTournier:2016:1}.

In this paper we:
\begin{itemize}
	\item derive the block version of \gcrodr, as implemented in \cite{PS.GCRODR-Code.2011} (\cf \Cref{section.deriv-impl}),
	\item discuss implementation decisions to favor block operations, including increasing the Krylov subspace block size
		using random vectors (\cf \Cref{subsection.impl}),
	\item demonstrate performance gains of a block Krylov subspace
		recycling over its single-vector counterpart (\cf
		\Cref{section.small-experiments,section.parameter-study}),
	\item show that the block, sparse matrix operations perform well in terms of appropriate data metrics relevant in the
		high-performance computing context,
	\item and study the efficiency of block operations specific to the block \gcrodr setting (\cf \Cref{section.data-movement}).
\end{itemize}
The last two goals are achieved through direct measurement of data movement and cache use efficiency on
the processor. Through carefully designed experiments, we are able to show that the cost of applying an operator to a block of
vectors is often marginally greater than applying the operator to a single vector in terms of data movement and usage cost
metrics.  Thus, we show that block methods can offer an accelerated convergence rate while reducing the overall data transmission
costs by avoiding data movement bottlenecks in modern hardware architectures.

\section{Preliminaries}\label{section.prelim}
When not necessary for the explanation or derivation of methods, we drop the index $i$ and consider the linear system
$\bA\paren[auto]{(}{)}{\bX_0 + \bT}=\bB$. Krylov subspace methods begin with a matrix $\bA\in\C^{n\times n}$ and a vector
$\bu\in\C^n$ and build a basis for the Krylov subspace
\begin{align}
	\cK_m(\bA,\bu)
	=
	\mathrm{span}\curly{\bu, \bA\bu,\ldots, \bA^{m-1}\bu}.
	\label{eqn.Krylov-basis}
\end{align}
We focus on methods that build an orthonormal basis for the Krylov subspace using
the Arnoldi process.  Let $\bV_m\in\C^{n\times m}$ be the
matrix with orthonormal columns generated by the Arnoldi process spanning $\cK_m(\bA,\bu)$.
Then we have the Arnoldi relation
\begin{align}
	\bA\bV_m
	=
	\bV_{m+1}\underline{\bH_m}
	\label{eqn.arnoldi-relation}
\end{align}
with $\underline{\bH_m}\in\C^{(m+1)\times m}$ upper Hessenberg; see, e.g.,
\cite[Section 6.3]{Saad.Iter.Meth.Sparse.2003} and
\cite{szyld.simoncini.survey.2007}.  

In general, for any $p\geq 1$, let $\bR_0=\bB-\bA\bX_0$ denote the $ n \times p$ initial residual. For the case $p=1$, at iteration $m$, we compute
$\bX_m = \bX_0 + \bT_m$, where $\bT_m\in\cK_m(\bA, \bR_0)$. In \gmres \cite{Saad.GMRES.1986}, we choose
\begin{align*}
	\bT_m
	=
	\argmin_{\bT\in \cK_m(\bA,\bR_0)}\norm[auto]{\bB-\bA(\bX_0 + \bT)}_2,
\end{align*}
and this is equivalent to solving the smaller minimization problem
\begin{align}
	\bY_m
	=
	\argmin_{\bY\in\C^{j}}\norm[auto]{ \underline{\bH_m}\bY - \norm[auto]{\bR_0}\be_1^{(m+1)}}_2,
	\label{eqn.GMRES-least-squares}
\end{align}
where we use the notation $\be_{\ell}^{(k)}$ to denote the $\ell$th Cartesian basis vector in $\R^{k}$,
and setting $\bX_m = \bX_0 + \bV_m\by_m$.  We call $\bT_m$ a \emph{correction}.
In restarted \gmres, i.e., (\gmres($m$)), we halt this process
at step $m$, discard the matrix $\bV_m$, and restart with the new initial residual
$\bR_0~\leftarrow~\bB - \bA \bX_m$.
This process is repeated until we achieve convergence.

\subsection{Block Krylov subspace methods}\label{subsection.block-krylov}
The extension of Krylov subspaces and the associated iterative methods to the block Krylov setting has been previously described
in,~e.g.,~\cite{Gutknecht2007,OLeary1980,Vital1990}. Though originally described for solving \eqref{eqn.Axb} in the case $p>1$,
such methods have also been proposed for accelerating convergence in the case that $p=1$. 

A block Krylov subspace $\K_m(\bA,\bR_0)$ is a generalization of the definition
of a Krylov subspace with more than one starting vector, i.e.,
\begin{align*}
	\K_m(\bA,\bR_0) 
	=
	\mathrm{blspan}\curly{\bR_0, \bA\bR_0, \bA^{2}\bR_0,\ldots, \bA^{m-1}\bR_0},
\end{align*}
where we note that by \emph{blspan}, we mean that we treat $n\times p$ vectors as a one-sided vector space, meaning that the linear
combinations are constructed using right-multiplication by $p\times p$ matrices.  Thus, elements of $
\K_m(\bA,\bR_0)$ have the form $\sum_{i=0}^{m-1}\bA^i\bR_0\bS_i$ where $\bS_i\in\C^{p\times p}$.  This one-sided vector space
approach is useful for understanding the behavior of these methods since it allows one to maintain the block structure in the
analysis by, e.g., considering the iteration in a space over the $^\ast$-algebra of $p\times p$ matrices.  This has been used to
great effect in, e.g., \cite{KubinovaSoodhalter:2020:1} which builds on ideas from \cite{Simoncini.Conv-Block-GMRES}.

It is straightforward to show that this interpretation is equivalent to treating this space as a vector space over $\C$ by observing
that
\begin{align}
	\K_m(\bA,\bR_0) 
	=&
	\cK_m(\bA,\br_0^{(1)}) +  \cK_m(\bA,\br_0^{(2)}) + \cdots  + \cK_m(\bA,\br_0^{(p)}),
	\nonumber
	\\
	=&
	\colspan\curly{\bR_0, \bA\bR_0, \bA^{2}\bR_0,\ldots, \bA^{m-1}\bR_0},
	\label{eqn.block-Krylov-sum-space}
\end{align}
where $\br_0^{(i)} = \bR_0\be_i$ is the $i$th column of $\bR_0$.  We mean this in the sense that for any vector $\bc\in\C^p$, we
can express $\paren[auto]{(}{)}{\sum_{i=0}^{m-1}\bA^i\bR_0\bS_i}\bc$ as a linear combination of elements from the constituent Krylov subspaces in
the sum \eqref{eqn.block-Krylov-sum-space}.

We denote by $L$ the block size used to generate the block Krylov subspace.  We consider two block Krylov subspace use-cases:  
\begin{itemize}
	\item if $p>1$, and we set $L=p$, we build a block Krylov subspace using $\bR_0$;
	\item if $p=1$ and $L>p$, we build the block Krylov subspace using $\bR_0\in\C^{n}$ and $L-1$ other
		vectors, independent from the residual.
\end{itemize}

Following the description in \cite[Section 6.12]{Saad.Iter.Meth.Sparse.2003}, we represent $\K_m(\bA,\bR_0)$ in terms of the block
Arnoldi basis $\curly{\bV_1,\bV_2,\ldots, \bV_m}$ where $\bV_i\in\C^{n\times p}$ has orthonormal columns and each column of
$\bV_i$ is orthogonal to all columns of $\bV_j$ for all $j\neq i$.  We obtain $\bV_1$ via the reduced QR-factorization
$\bR_0~=~\bV_1\bF_0$ where $\bF_0\in\C^{p\times p}$ is upper triangular.  We can generate $\bV_{m+1}$ with the
block Arnoldi step; see, \eg, \cite[Algorithm 6.22]{Saad.Iter.Meth.Sparse.2003}.  Let $\bW_m = \begin{bmatrix}\bV_1 & \bV_2 \cdots
\bV_m\end{bmatrix}\in\C^{n\times mp}$. Let $\underline{\bH_m} = (\bH_{ij})\in\C^{(m+1)p\times mp}$. This yields the block Arnoldi
relation
\begin{align}
	\bA\bW_m 
	=
	\bW_{m+1}\underline{\bH_m},
	\label{eqn.block-arnoldi-relation}
\end{align}
where $\underline{\bH_m}\in\C^{(m+1)p\times mp}$ is block upper Hessenberg, with $p\times p$ blocks,
in which the lower block-subdiagonal is composed of upper-triangular matrices.
A straightforward generalization of \gmres for block Krylov subspaces (called block
\gmres), first described in \cite{Vital1990}; see, e.g., \cite[Chapter 6]{Saad.Iter.Meth.Sparse.2003} for
more details.  For $p>1$, one solves the generalization
of the single-vector \gmres minimization problem,
\begin{align}
	\bY_m 
	=
	\argmin_{\bY\in\C^{mp\times p}}\norm[auto]{\underline{\bH_m}\bY - \bE_1^{(m+1)p}\bF_0}_F
	\label{eqn.true-block-gmres-min}
\end{align}
and setting $\bX_m = \bX_0 + \bW_m\bY_m$ where $\bE_1^{(m+1)p}$ is
the matrix containing the first $p$ columns of the order $(m+1)p$ identity matrix.  It is easy to show this is equivalent to
computing the minimum residual $2$-norm correction over the block Krylov subspace, one column at-a-time.

In the case that $p=1$ with the block size having been enlarged to $L$, the subspace is $\K_m(\bA,\bV_1)$ where $\bR_0 =
\norm[auto]{\bR_0}\bV_1\be_1^{(L)}$.  At iteration $m$, one solves
\begin{align}
	\bY_m 
	=
	\argmin_{\bY\in\C^{mL}}\norm[auto]{\underline{\bH_m}\bY -  \bE_1^{(m+1)L}\bF_0\be_1^{(L)} }_2
	\label{eqn.1rhs-block-gmres-min}
\end{align}
where $\be_1^{(L)}\in\curly{0,1}^{L}$ is the first column of the identity matrix. One then sets $\bX_m= \bX_0 + \bW_m\bY_m$ as before. In
this case, one minimizes the single-vector residual over the block Krylov subspace.

The core message is that block \gmres is a residual minimization method generalizing \gmres to the block Krylov subspace setting.
By understanding recycled \gmres as being from a class of methods that minimize over the sum of two subspaces (one of which
is a Krylov subspace), we are able to chart a clear path forward for extending \gcro-based recycled \gmres to the block Krylov subspace
setting.

\subsection{Recycled \gmres}\label{section.rgmres}
\emph{Subspace recycling} is a type of augmented Krylov subspace method wherein one augments a Krylov subspace with vectors
generated by a previous iteration for the same or a previous system,\footnote{or from some other
helpful source}  a technique first denoted as \emph{recycling} in \cite{Parks.deSturler.GCRODR.2005}.
An augmented Krylov
subspace method is such that  a correction to the initial approximation is
computed not just over a Krylov subspace $\cK$ but instead over
an \emph{augmented} Krylov subspace of the form $\cU + \cK$, where $\cU$ is available before the start of the iteration.  We use
the name recycled \gmres (\rgmres) to encompass all augmented Krylov subspace methods that minimize the residual norm over
an augmented Krylov subspace.  The most successful implementations are those presented in the \gcro framework introduced in
\cite{deSturler.GCRO.1996}. It is shown that such methods equivalently can be expressed as a \gmres iteration applied to the
original linear system, left-multiplied with a specially chosen projector.  The survey \cite{SoodhalterDeSturlerKilmer:2020:1}
goes into much more detail and generalizes the idea to augmented/recycling approaches not based on the minimization of an error
functional. 
What differentiates the methods within this class is how the augmenting subspace
$\cU$ is computed and updated.

We briefly review the  method described in \cite{Parks.deSturler.GCRODR.2005}. For simplicity, we continue to drop the superscript
$^{[i]}$.  We assume there is an augmentation space $\cU$ that is available before the start of the iteration. This algorithm
represents the combination of two approaches: those originating from the implicitly restarted Arnoldi method \cite{LS.1996}, such
as Morgan's \gmresdr \cite{Morgan.GMRESDR.2002}, and those descending from de Sturler's \gcro method \cite{deSturler.GCRO.1996}.
\gmresdr is a restarted \gmres-type algorithm, where at the end of each cycle, harmonic Ritz vectors are computed, and a subset of
them are used to augment the Krylov subspace generated at the next cycle. The \gcro method allows the user to select the optimal
correction over arbitrary subspaces. This concept is extended by de Sturler in~\cite{deSturler.GCROT.1999} (and simplified in
\cite{Manteuffel.LGMRES.2005}), where a framework is provided for selecting the optimal subspace to retain from one cycle to the
next so as to minimize the error produced by discarding useful information accumulated in the subspace for candidate solutions
before restart. This algorithm is called \gcrot, and this procedure is referred to as \quotes{optimal truncation}. Parks et~al.\
in \cite{Parks.deSturler.GCRODR.2005} combine the ideas of  \cite{Morgan.GMRESDR.2002} and \cite{deSturler.GCROT.1999} and extend
them to a sequence of slowly-changing linear systems and recycling with harmonic Ritz vectors. They call their method \gcrodr.

All methods that minimize the residual over an augmented subspace of the form $\cU + \fK_m$ (i.e., a fixed space $\cU$ and an
iteratively generated space $\fK_m$ that increases in dimension at each iteration, with $d_m:=\dim\fK_m$) have common structural characteristics that are
exploited when designing an \rgmres algorithm.  \gcro-based approaches such as \gcrodr are demonstrations of this point. We
distill the most important aspects of the general theory presented in the survey \cite{SoodhalterDeSturlerKilmer:2020:1} in the
residual minimization setting.

Consider solving $\bA\paren[auto]{(}{)}{\bX_0 + \bS + \bT} = \bB$ with $\bR_0=\bB-\bA\bX_0$ for $p\geq 1$.  In \cite{SoodhalterDeSturlerKilmer:2020:1} 
the authors express the augmented iterative method as constructing the approximation $\bS+\bT\approx \bS_m + \bT_m\in\cU +
\fK_m$.  The augmented residual minimization approach selects $\bS_m\in\cU$ and $\bT_m\in\fK_m$ such that $\norm[auto]{\bB -
\bA\paren[auto]{(}{)}{\bX_0 + \bS_m + \bT_m}}$ is minimized. Let us express $\bS_m = \bU \bZ_m$ where $\bU\in\C^{n\times k}$ has
columns spanning $\cU$ and
$\bZ_m\in\C^{k\times p}$, and $\bT_m=\bW_m\bY_m$ where $\bW_m\in\C^{n\times d_m}$ has columns spanning $\fK_m$ and
$\bY_m\in\C^{d_m\times p}$.  A key observation from \cite[Section 5.1.1]{SoodhalterDeSturlerKilmer:2020:1} is that solving 
\begin{align*}
	\bA
	\begin{bmatrix}
		\bU & \bW_m
	\end{bmatrix}
	\begin{bmatrix}
		\bZ_m
		\\
		\bY_m
	\end{bmatrix}
	\approx
	\bR_0
\end{align*}
via a least-squares approach using the normal equations is equivalent to approximating the solution of the singular, consistent
linear system
\begin{align}
	\paren[auto]{(}{)}{\bI - \bPhi}\bA\bT
	=
	\paren[auto]{(}{)}{\bI - \bPhi}
	\bR_0
	\label{eqn.projected-subproblem}
\end{align}
with $\bT_m\in\fK_m$ via residual minimization and then constructing 
\begin{align}
	\bX_m 
	= 
	\bX_0 
	- 
	\bPi\bT 
	+ 
	\paren[auto]{(}{)}{\bI - \bPi}\bT_m, 
	\label{eqn.full-rgmres-solution}
\end{align}
where $\bPhi$ is the orthogonal projector onto $\cC := \bA\cU$, $\bPi$ is the $\bA^\ast\bA$-orthogonal projector onto
$\cU$, and the resulting initial error projection $\bPi\bT$ is a computable quantity.  Furthermore, it was shown that for any such
augmented subspace residual minimization, the full residual and the residual for the projected subproblem are the same, i.e.,
\begin{align}
	\bB - \bA\bX_m
	=
	\paren[auto]{(}{)}{\bI - \bPhi}\bR_0 - \paren[auto]{(}{)}{\bI - \bPhi}\bA\bT_m.
	\label{eqn.full-projected-residual-equivalent}
\end{align}
\begin{remark}
	We note that the framework developed in \cite{SoodhalterDeSturlerKilmer:2020:1} is independent of the choices of
	subspaces $\cU$ and $\fK_m$.  The projected subproblem \eqref{eqn.projected-subproblem} comes from performing a Galerkin
	residual minimization over a sum of two subspaces.  
\end{remark}
Indeed, it is the choice of $\fK_m$ that takes a general augmented subspace residual minimization and turns it into a \gcro-based
\rgmres algorithm.  If we specify $\fK_m = \cK_m\paren[auto]{(}{)}{\paren[auto]{(}{)}{\bI - \bPhi}\bA, \paren[auto]{(}{)}{\bI -
\bPhi}\bR_0}$, then this process becomes equivalent to applying a GMRES iteration directly to \eqref{eqn.projected-subproblem}. It
follows then from \eqref{eqn.full-projected-residual-equivalent} that the residual convergence behavior of a \gcro-based \rgmres
method is governed completely by the behavior of \gmres applied to the projected subproblem \eqref{eqn.projected-subproblem}; \cf
\Cref{section.conv} for further discussion of this fact.

\subsubsection{Standard \gcrodr version of \rgmres} Suppose we are solving \eqref{eqn.Axb} with $p=1$, and we have a
$k$-dimensional subspace $\cU$ which is spanned by vectors recycled either from a previous linear system solve or in the previous
iteration cycle and whose image under the action of $\bA$ is $\cC = \bA\,\cU$.  Let $\bPhi$ be the orthogonal projector onto $\cC$.
As discussed in \Cref{section.rgmres}, minimizing the residual over an augmented Krylov subspace $\cU + \fK_m$ is equivalent to
applying a \gmres to~\eqref{eqn.projected-subproblem}.
We generate the Krylov subspace using a projected version of the Arnoldi
process.  After $m$ iterations, \gmres applied to \eqref{eqn.projected-subproblem} produces the correction
$\bT_m\in\cK_m\paren[auto]{(}{)}{\paren[auto]{(}{)}{\bI - \bPhi}\bA,\paren[auto]{(}{)}{\bI - \bPhi}\bR_0}$. At the end of the cycle, an updated $\cU$ is constructed,
the Krylov subspace basis is discarded, and we restart. At convergence, $\cU$ is saved, to be used when solving the next linear
system.

Practical construction of \eqref{eqn.full-rgmres-solution} is straightforward. A fundamental choice for implementing any recycling
method is the choice of bases for $\cU$ and $\cC$.  For \gcro-based \rgmres implementations, the usual choice is to maintain an
orthonormal basis for $\cC$.  This greatly simplifies the representation of the projectors $\bPi$ and $\bPhi$ since
$\bC^\ast\bC=\bI$.  It follows that  
\begin{align*}	
	\bPi 
	= 
	\bU\paren[auto]{(}{)}{\bU^\ast\bA^\ast\bA\bU}^{-1}\bU^\ast\bA^\ast\bA
	=
	\bU\bU^\ast\bA^\ast\bA
	=
	\bU\bC^\ast\bA,
	\qquad
	\mbox{and}
	\qquad
	\bPhi
	=
	\bC\bC^\ast.
\end{align*}
It is well documented that all recycling methods that fit into the framework described in \cite{SoodhalterDeSturlerKilmer:2020:1} have
a projected error term of the form $\bPi\bT$ that is practically computable.  In the \rgmres setting, this follows from the
structure of $\bPi$ since $\bPi\bT =  \bU\bC^\ast\bR_0$.  When applying \gmres to \eqref{eqn.projected-subproblem}, we obtain the
modified (projected) Arnoldi relation
\begin{align}
	\paren[auto]{(}{)}{\bI - \bPhi}\bA\bV_m 
	=& 
	\bV_{m+1}\underline{\bH_m}
	\\
	\iff
	\bA\bV_m
	=&
	\bPhi\bA\bV_m
	+
	\bV_{m+1}\underline{\bH_m}
	=
	\bC\bB_m
	+
	\bV_{m+1}\underline{\bH_m},
	\label{eqn.modified-arnoldi}
\end{align}
where $\bB_m = \bC^\ast\bA\bV_m$.  
The action of $\paren[auto]{(}{)}{\bI - \bPhi}$ is implemented as
an orthogonalization away from the orthonormal columns of $\bC$, with the coefficients stored in $\bB_m$.  Solving the usual
\gmres minimization, we obtain $\bY_m = \argmin_{\bY\in\C^m}\norm[auto]{\underline{\bH_m}\bY - \beta\be_1}$ where $\beta =
\norm[auto]{\paren[auto]{(}{)}{\bI - \bPhi}\bR_0}$, and $\bT_m = \bV_m\bY_m$.  Lastly, we obtain $\bPi\bT_m$ by observing that
\begin{align*}
	\bPi\bT_m
	=
	\bU\bC^\ast\bA\bV_m\bY_m
	=
	\bU\bB_m\bY_m,
\end{align*}
which involves already-computed quantities.

Convergence analysis for augmented Krylov subspace methods was previously presented in, e.g.,
\cite{Eiermann.Analysis-accel.2000, Saad.Deflated-Aug-Krylov.1997}. In the context of \rgmres, the thesis of Gaul
\cite{Gaul.2014-phd} and the references therein are all excellent sources on this topic. 

Iterating orthogonally to an approximate invariant subspace to accelerate convergence of \gmres can be justified by the
theoretical work in \cite{Simoncini2005}, wherein it is shown that the widely observed superlinear convergence behavior of \gmres,
is governed by how well the Krylov subspace approximates a certain
invariant subspace of $\bA$.
This analysis complements previous discussions of this superlinear convergence phenomenon; 
see e.g.,~\cite{Ipsen.GMRES-minimal-poly.1996,Vorst1993a}.  
However, it should be noted that the theory describing the effectiveness of \rgmres
applied to a non-normal system is not yet fully understood.  Explanations that characterize acceleration of convergence in terms of
projections onto invariant subspaces do not take into consideration the non-normality of the matrix; and, thus, the effects of
ill-conditioning of the eigenbasis are ignored.  This has been mentioned in  \cite{desturler.recycling-convergence-abstract.2012},
but has not yet been fully explored.

\section{Recycled Block \gmres}\label{section.deriv-impl}

The framework discussed in \Cref{section.rgmres} is compatible with any minimum residual iterative
method over an augmented Krylov subspace.  We simply let $\fK_m$ be a block Krylov subspace.  We describe the method to accommodate
$p\geq 1$ and describe the differences when working with a true block method ($p>1$).
In that vein,
we consider a generic block size $L$.

Given a subspace $\cU$ we derive the block recycled \gmres iteration thusly. Using the block Arnoldi process, we generate a basis
for the subspace \linebreak $\fK_m=\K_m\paren[auto]{(}{)}{\paren[auto]{(}{)}{\bI - \bPhi}\bA, \paren[auto]{(}{)}{\bI - \bPhi}\bR_0}$ where
the orthonormal columns of $\bW_m\in\C^{n\times mL}$ span the subspace. By construction, the columns of $\bW_{m+1}$ are orthogonal
to the columns of $\bC$, yielding a block version of \eqref{eqn.modified-arnoldi},
\begin{align}
	\bA\bW_m
	=
	\bC\bB_m
	+
	\bW_{m+1}\underline{\bH_m}
	\label{eqn.block-modified-arnoldi}
\end{align}
where $\bB_m = \bC^{\ast}\bA\bW_m\in\C^{k\times mL}$ represents the entries generated by orthogonalizing the columns of the new
block Krylov basis vector against $\cC$.  The derivation proceeds just as in \Cref{section.rgmres}.  It still holds that
$\bPi\bT=\bU\bC^\ast\bR_0$.  We obtain $\bT_m$ as the $m$th block \gmres approximation to the solution of
\eqref{eqn.projected-subproblem}, and $\bPi\bT_m = \bC\bB_m\bY_m$.  We obtain $\bY_m$ from the minimization \eqref{eqn.true-block-gmres-min} in the case
that $p>1$ and \eqref{eqn.1rhs-block-gmres-min} in the case that $p=1$.

\begin{proposition}\label{prop.residual-equiv}
The block residual produced by block \gcrodr applied to \eqref{eqn.Axb} and that produced by block \gmres applied to
\eqref{eqn.projected-subproblem} are equal; i.e., 
	\begin{align*}
		\bB - \bA\paren[auto]{(}{)}{\bX_0 + \bS_m + \bT_m}
		=
		\paren[auto]{(}{)}{\bI - \bPhi}\paren[auto]{(}{)}{\bR_0 - \bA\bT_m}.
	\end{align*}
\end{proposition}
\begin{proof}
	At iteration $m$ of block \gcrodr, we have the block residual
\begin{align*}
	\bR_m
	=
	\bB - \bA\paren[auto]{(}{)}{\bX_0 + \bS_m + \bT_m}.
\end{align*}
Inserting the expressions from earlier for the two corrections, we can write
\begin{align*}
	\bR_m 
	=& 
	\bR_0 - \bA\paren[auto]{(}{)}{-\bU\bB_m\bY_m + \bW_m\bY_m}\\
    	=& 
	\bR_0 + \bC\bB_m\bY_m - \begin{bmatrix} \bC & \bW_{m+1} \end{bmatrix}\begin{bmatrix}\bB_m\\\underline{\bH_m}\end{bmatrix}\bY_m\\
    	=& 
	\bR_0 + \bC\bB_m\bY_m - \bC\bB_m\bY_m - \bW_{m+1}\underline{\bH_m}\bY_m
    	=
	\bR_0 - \bW_{m+1}\underline{\bH_m}\bY_m.
  \end{align*}
  This is the block residual produced by block \gmres applied to the projected problem \eqref{eqn.projected-subproblem},
  proving the proposition.
\end{proof}

One can, in fact, represent \eqref{eqn.block-modified-arnoldi} as one large blocked Hessenberg relation. 
Let
\begin{align}
	\widehat{\bW}_m
	=
	\begin{bmatrix}
		\bU & \bW_m
	\end{bmatrix},
	\qquad\widetilde{\bW}_{m+1}
	=
	\begin{bmatrix}
		\bC & \bW_{m+1}
	\end{bmatrix},
	\ \ 
	\mbox{and}
	\ \ 
	\underline{\bG_m}
	=
	\begin{bmatrix}
		\bI_k & \bB_m 
		\\ 
		\b0 & \underline{\bH_m}
	\end{bmatrix},
	\label{eqn.mod-arnoldi-notation}
\end{align}
with $\bI_k$ being the $k\times k$ identity matrix. It follows that we can write
\begin{align}
	\bA \widehat{\bW}_m = \widetilde{\bW}_{m+1}\underline{\bG_m}.
	\label{eqn.mod-arn-rel}
\end{align}
We do not advocate implementing a \gcrodr method using a compact augmented Arnoldi relation, as it complicates the
algorithm and introduces possible stability issues; see, e.g., \cite{Parks.deSturler.GCRODR.2005}. However, it is useful to
introduce it for the computation of harmonic Ritz vectors, \cf \Cref{section.ritz-comp}.
\begin{algorithm}[tb!]
\caption{The Block \rgmres Algorithm} 
\label{alg.block-GCRODR}
\SetKwInOut{Input}{Input}\SetKwInOut{Output}{Output}
\Input{$\bA\in \C^{n\times n}$, $\bB\in\C^{n\times p}$, $\bX_0\in\C^{n\times p}$, $\varepsilon > 0$ the convergence tolerance, $m$ the number of block Arnoldi vectors generated, $k$ the desired dimension of the space to be recycled,
and possibly $\bU$, whose $k$ columns span the space to be recycled (optional, if available).}
\Output{$\bX\in\C^{n\times p}$ such that $\norm[auto]{\bB - \bA\bX}_F\leq \varepsilon$}
Set $\bR_0 = \bB - \bA\bX_0$, Set $i = 0$\\
Set $\bE = \begin{bmatrix}\be_1^{(mp)} & \cdots & \be_p^{(mp)}\end{bmatrix}$\\
\If{$\bU$ is available (e.g., defined from solving a previous linear system)}{
	Define $\bC$ using reduced QR Factorization $\bA\bU=\bC\bS$, $\bU\leftarrow \bU\bS^{-1}$, and $\bPhi = \bC\bC^\ast$\\
	$\bX_0 \leftarrow \bX_0 + \bU\bC^{\ast}\bR_0$, $\bR_0 \leftarrow \bR_0 - \bPhi\bR_0$
}\Else{
	Define $\bV_1$ using reduced QR Factorization $\bR_0 = \bV_1\bZ$\\
	Perform $m$ steps of block \gmres, generating~$\bW_{m+1}$~and~$\underline{\bH_m}$ and solve 
${\bY}_0 = \argmin_{\bY \in \C^{mp\times p}}\norm[auto]{\underline{\bH_m}\bY - \bE\bZ }_F$, \\
	$\bX_0 \leftarrow \bX_0 + \bW_m{\bY}_0$, 
$\bR_0 \leftarrow \bR_0 - \bW_{m+1}\underline{\bH_m}{\bY}_0$\\
	Select a subspace $\cU$ of $\cK_m(\bA,\bV_1)$ to recycle. \\
	Compute $\bU$ having basis vectors of $\cU$ as columns such that $\bC=\bA\bU$
	has orthonormal columns, and $\bPhi = \bC\bC^\ast$
}
\While{$\norm[auto]{\bB - \bA\bX_i}_F>\varepsilon$}{
	$i\leftarrow i+1$\\
	Define $\bV_1$ using reduced QR Factorization $\bR_{i-1} = \bV_1\bZ$\\
	Compute a basis for $\K_m((\bI - \bPhi)\bA,\bV_1)$ using the block Arnoldi method, generating $\bW_{m+1}$, $\underline{\bH_m}$, and $\bF_m$.\\
	Solve the block \gmres least-squares subproblem $\bY_i = \argmin_{\bY\in\C^{mp\times p}}
\norm[auto]{\underline{\bH_m}\bY - \bE\bZ }_F$\\
	Set $\bX_i = \bX_{i-1} + \bW_m\bY_i - \bU\bF_m\bY_i$, Set $\bR_i = \bR_{i-1} - \bW_{m+1}\underline{\bH_m}\bY_i$\\
	Define $\bD$ to be the diagonal matrix such that $\widetilde{\bU} = \bU\bD$ has columns of unit norm.\\
	Set $\underline{\bG_m} = \begin{bmatrix}\bD & \bF_m \\ \b0 & \underline{\bH_m}\end{bmatrix}$\tcp*{Scaling $\bU$ with $\bD$ for stability}
	Compute $\cU_{new} \subset\cU + \K_m((\bI-\bPhi)\bA,\bV_1)$, $\cU \leftarrow \cU_{new}$\\
  Compute $\bU$ such that $\bC\leftarrow\bA\bU$
	has orthonormal columns, and $\bPhi = \bC\bC^\ast$
}
Store $\bU$ in memory to serve as initial recycle subspace for next function call.
\end{algorithm}

We use $\underline{\bG_m}$ to compute a new approximate invariant subspace; and if we have not converged,
we begin the next cycle.  Algorithm \ref{alg.block-GCRODR} gives a complete pseudocode description of the
algorithm.

\subsection{Implementation considerations}\label{subsection.impl}
We discuss implementation decisions made in light of the fact that the proposed method is built upon a block Krylov subspace
method.

\subsubsection{Householder reflection storage}
Working with a block Hessenberg matrix introduces some additional computational challenges as compared to the non-block case.  We
elaborate on our approach to block triangularization of the block Hessenberg matrix. In the case of $L=p=1$ the upper Hessenberg
matrix $\underline{\bH_m}$ has only one subdiagonal entry per column.  To compute its QR-factorization at each step of the method,
one annihilates the subdiagonal entry of each column in a progressive manner using Givens rotations, which are retained compactly
in the form of sines and cosines to be applied to subsequent columns. For $L>1$, $\bH_m$ is block upper Hessenberg.  For a
block of columns,  newly generated by a step of the block Arnoldi procedure, new Householder reflections are computed
column-by-column.  However, we must first apply all previously generated reflections to this new block of columns.  We employ the
strategy of Gutknecht and Schmelzer \cite{Gutknecht2008}. One stores the Householder reflections for a block column as a single
matrix and applies them all at once.  This exchanges $p$ applications of previous Householder reflections for one dense
matrix-matrix multiplication.  This dense matrix-matrix multiplication can be performed as a level-3 \blas operation.  It has been
noted that for certain approaches to understanding the behavior of block \gmres, the block Householder transformations can be
difficult to interpret.  Indeed, in \cite{KubinovaSoodhalter:2020:1}, 
the block GMRES iteration for the case
$p>1$ is interpreted in terms of a block vector iteration over a one-sided vector space with scalars from the $^\ast$-algebra of $\C^{p\times
p}$, and this interpretation is used to meaningfully extend the results from \cite{Greenbaum.Any-Curve-Possible-GMRES.1996} to
the block \gmres setting.  When considering a block Krylov iteration in this way, it is more natural to formally
consider\footnote{We say only formally because it is not practical to implement them.} the triangularization of $\underline{\bH_m}$ via a
block generalization of Givens rotations.

For a block version of \gcrodr, the computation and updating of the recycled subspace $\cU$ is a direct generalization of the
non-block case. We note that if no space $\cU$ is given at execution, we follow \cite{Parks.deSturler.GCRODR.2005} and run a cycle
of block \gmres, computing harmonic Ritz vectors with respect to the block Krylov subspace at the end of the cycle. We discuss next
this computation in more detail and also ponder other recycling strategies.

\subsubsection{Harmonic Ritz vector computation}\label{section.ritz-comp} 
This is the
strategy implemented in \cite{PS.GCRODR-Code.2011} following the harmonic Ritz vector deflation strategy 
in
\cite{Morgan.GMRESDR.2002}. At the end of the cycle, we generated an orthonormal basis for the subspace $\K_m(\bA,\bR_0)$
with the block Arnoldi relation $\bA\bW_m = \bW_{m+1}\underline{\bH_m}$.  Following \cite{Morgan.Restarted-GMRES-eig.1995}, the
block harmonic Ritz problem for $\K_m(\bA,\bR_0)$ is to find all pairs 
\begin{align}
	(\by,\mu)
	\in
	\K_m(\bA,\bR_0)
	\times
	\C
	\mbox{ such that }
	\bA^{-1}\by - \mu\by \perp \bA\K_m(\bA,\bR_0).
	\label{eqn.block-harm-orth}
\end{align}
As with the scalar case, \eqref{eqn.block-harm-orth} can be equivalently solved as a generalized eigenvalue problem whose solution
pairs $(\bt,\tilde{\theta})\in\C^{mp}\times\C$ can be used to reconstruct the pairs $(\by,\mu)\in\C^{n}\times \C$ as described in
the following.
\begin{proposition}\label{prop.block-harm-ritz} 
	Given the block Krylov subspace $\K_m(\bA,\bR_0)$, solving the harmonic Ritz problem \eqref{eqn.block-harm-orth} is
	equivalent to solving the $mp \times mp$ eigenvalue problem
\begin{align*}
	(\bH_m + (\bH_m^{\ast})^{-1}\widehat{\bE}(\bH_{m+1,m}^{\ast}\bH_{m+1,m})\widehat{\bE}^{\ast})\bt
	=
	\tilde{\theta}\bt
\end{align*}
 and then for a solution pair $(\bt,\tilde{\theta})$ assigning $\by = \bW_m\bt$, where the columns of
	$\widehat{\bE}\in\curly{0,1}^{mp\times p}$ are columns $(mp-p+1),\ldots , mp$ of the identity matrix of order $mp$, and
	$\mu = 1/\tilde{\theta}$.
\end{proposition}

\noindent It should be noted that, as a practical matter, the expression
\begin{align*}
	\bH_m + (\bH_m^{\ast})^{-1}\widehat{\bE}(\bH_{m+1,m}^{\ast}\bH_{m+1,m})\widehat{\bE}^{\ast}
\end{align*}
simply means that the last $p$ columns of $\bH_m$ are modified by the $mp\times p$ matrix
\begin{align*}
	(\bH_m^{\ast})^{-1}\widehat{\bE}(\bH_{m+1,m}^{\ast}\bH_{m+1,m}).
\end{align*}
\begin{proof}[Proof of Proposition~\ref{prop.block-harm-ritz}]
This is a generalization of the harmonic Ritz computation in the case of a single-vector Krylov subspace; see e.g., \cite{Morgan2000}.
We can prove this through algebraic manipulation using the block Arnoldi relation.  Condition \eqref{eqn.block-harm-orth} is equivalent to
\begin{align*}
	(\bA\bW_m)^{\ast}(\bA^{-1}\by - \mu\by) 
	=&
	0\\
	(\bW_{m+1}\underline{\bH_m})^{\ast}(\bA^{-1}\bA\bW_m\bt - \mu\bA\bW_m\bt)
	=&
	0  \\
	\underline{\bH_m}^{\ast}\bW_{m+1}^{\ast}(\bW_m\bt - \mu\bW_{m+1}\underline{\bH_m}\bt)
	=&
	0  \\
	\underline{\bH_m}^{\ast}\bW_{m+1}^{\ast}\bW_m\bt
	=&
	\mu \underline{\bH_m}^{\ast}\underline{\bH_m}\bt \\
	\bH_m^{\ast}\bt
	=&
	\mu (\bH_m^{\ast}\bH_m + \widehat{\bE}(\bH_{m+1,m}^{\ast}\bH_{m+1,m})\widehat{\bE}^{\ast})\bt \\
	\tilde{\theta}\bt
	=&
	(\bH_m + (\bH_m^{\ast})^{-1}\widehat{\bE}(\bH_{m+1,m}^{\ast}\bH_{m+1,m})\widehat{\bE}^{\ast})\bt
\end{align*}
\end{proof}

The computation in the case of the augmented subspace $\cU + \cK_m(\paren[auto]{(}{)}{\bI - \bPhi}\bA,\paren[auto]{(}{)}{\bI - \bPhi}\bR_0)$ is similar to the computation employed in
\cite{Parks.deSturler.GCRODR.2005}, as described in the following result, whose proof
is nearly identical to the one developed in \cite[Equation 2.16]{Parks.deSturler.GCRODR.2005}. 
\begin{proposition} 
In a cycle of block recycled \gmres, if we have generated an augmented
	space $\cU + \cK_m(\paren[auto]{(}{)}{\bI - \bPhi}\bA,\paren[auto]{(}{)}{\bI - \bPhi}\bR_0)$ then solving the associated harmonic Ritz problem is equivalent to solving the
	generalized eigenvalue problem
\begin{align}
	\underline{\bG_m}^{\ast}\underline{\bG_m}\bt 
	= 
	\tilde{\theta}\underline{\bG_m}^{\ast}\widetilde{\bW}_{m+1}^{\ast}\widehat{\bW}_m\bt
	\label{eqn.harm-ritz-gen-eig}
\end{align}
and assigning $\by = \widehat{\bW}_m\bt$ for each solution pair $(\tilde{\theta},\bt)$
where $\underline{\bG_m}$, $\widetilde{\bW}_{m+1}^{\ast}$, and $\widehat{\bW}_m$ are
defined as in \eqref{eqn.mod-arnoldi-notation}.
\end{proposition}

\subsubsection{Other recycled space selection techniques}

Block Krylov subspaces have been originally proposed for the computation of eigenvalues/eigenvectors with the justification that
they generate richer subspaces, see, e.g., \cite{RMNW.2010,Ruhe1979}.  Thus, recycling approximate eigenvectors, as in
\cite{RMNW.2010,A-Giraud-YanFei.GM-Inexact.2014,Darnell2008,Morgan.GMRESDR.2002,Parks.deSturler.GCRODR.2005} offers the
possibility of rapidly acquiring high quality eigenvector approximations with which to deflate.  This makes the use of block
\gcrodr or some block/non-block hybrid strategy more attractive for the case $p=1,\ L>1$.  For the first few cycles, one can
inflate the block size in order to more quickly obtain a high quality recycled subspace $\cU$ and then switch at some restart to
non-block \gcrodr (i.e., $L=1$) thereafter.

However, our motivation arises mainly from considerations
in the high-performance computing setting.  For dense linear algebra computations, it has been shown that level-3 \blas (i.e.,
matrix-times-matrix) operations exhibit superior data movement efficiency properties, as measured amount of data moved per
operation and efficiency of data reuse in cache \cite{GH.Ng.OP.Plemmons.RSV.1990-book}.  The assumption in designing this
algorithm is that sparse matrix-times-matrix operations would also exhibit similar superior properties and that level-1 and
level-2 \blas operations generalize to level-3 \blas. This has been previously
discussed \cite{H.2010,MHDY.2009}. Careful experimentation will be necessary to demonstrate this, not only to understand this
behavior for the application of a large, sparse operator $\bA$ but also for the application of the projected operator $\prn{\bI -
\bPhi}\bA$.

In the current version of our codes \cite{PS.GCRODR-Code.2011}  (as well as in the current version of the publicly available \gcrodr codes \cite{Parks:2004:1})
harmonic Ritz vectors are computed to generate a subspace
to recycle.  

Indeed, there are other recycling strategies discussed in the literature. Morgan suggests that in an eigenvector deflation
algorithm based upon \fom, called \fomdr \cite{Morgan.GMRESDR.2002}, deflation using Ritz vectors is more effective.  It is
suggested in \cite{Parks.deSturler.GCRODR.2005} that perhaps a mix of Ritz and harmonic Ritz vectors may be appropriate in
some cases.  In his paper on optimal truncation methods \cite{deSturler.GCROT.1999}, de Sturler demonstrates that one can
calculate which subspace of dimension $k$ of the current Krylov subspace of dimension $m$ most important to maintain orthogonality
against, for the purpose of reducing the residual. This subspace is then recycled under the assumption that it is most important
to continue to maintain orthogonality with respect to this subspace.  Ahuja et al.\ \cite{Ahuja2010}, observed that the
preconditioned systems with which they dealt had eigenvalue clusters well separated from the origin, rendering the use of harmonic
Ritz vectors less effective.  Instead, they chose to recycle Krylov vectors which had dominant components in the right-hand side,
and this gave improved convergence results.

Gaul and Schl\"ommer \cite{GS.2013-arXiv} suggest that in the context of recycled \minres being used to solve a
Schr\"odinger-type equation, Ritz vectors are good candidates with which to recycle.  In \cite{CF.Tuminaro.Recyc-POD.2015}, the
authors propose a method of recycling using a proper orthogonal decomposition approach coming from model order reduction. In the
context of ill-posed image recovery problems, it has been demonstrated that one can also augment the Krylov subspace with vectors
which encode knowledge of characteristics of the true solution, e.g., edge characteristics of the image \cite{MRH.2014}.  In that
work, flexible \gmres \cite{Saad1993} is used to augment the subspace.  This follows from the work in
\cite{BR-2.2007,BR.2007,DGH.2014} in which \gmres for ill-posed problem is augmented with vectors encoding features of the
reconstructed image which are difficult for a Krylov method to reconstruct (such as discontinuities and hard edges). Using the
augmented method framework discussed in \cite{SoodhalterDeSturlerKilmer:2020:1}, the author of \cite{Soodhalter:2022:1} re-interpreted
the work of \cite{DGH.2014} in order to propose an alternative implementation.

\section{Convergence discussion}\label{section.conv}

It is shown, e.g., in \cite{Simoncini2005}, that the convergence of \gmres accelerates, entering a superlinear
phase, once the Krylov method has adequately captured a subspace spanned by eigenvectors associated to
eigenvalues which often cause slow convergence, i.e., those near the origin.
For these eigenvalues, low-degree residual polynomial interpolation can be difficult;
see also, \cite{Ipsen.GMRES-minimal-poly.1996,Vorst1993a}. This explains some of
the convergence difficulties exhibited by restarted methods, in which we discard the entire basis and start
over.  Furthermore, once this eigenspace is well-represented by the Krylov subspace,
\cf \cite{Simoncini.Nonoptimal-bases.2005}, the convergence behavior mimics that of an
operator from which the eigenspace has been removed.  This is one motivation for the subspace
augmentation and recycling technique, e.g., \cite{Morgan.GMRESDR.2002,Parks.deSturler.GCRODR.2005}.
By recycling a selected subspace and iterating orthogonally to it, we hope to enter the superlinear convergence
phase of \gmres earlier.  By building a block Krylov subspace, one can capture these invariant subspaces in
fewer iterations.

We mention that it is well understood that  the eigenvalues themselves can have no connection to the residual convergence
pattern of GMRES \cite{Greenbaum.Any-Curve-Possible-GMRES.1996}, a result that has been extended to the block \gmres setting
\cite{KubinovaSoodhalter:2020:1}. Examples are presented in \cite{CarsonLiesenStrakos:2024:1} that illuminate
the complicated nature of the mechanics of \gmres convergence speed. 

As we have shown, block \rgmres iteration is equivalent to a block \gmres iteration applied to a projected problem.
Thus, the convergence results for block \gmres can be extended to the recycled block \gmres case.
We focus without loss of generality on the true block method case of $p>1$.
Simoncini and Gallopoulos \cite{Simoncini.Conv-Block-GMRES}
discussed the convergence properties of block \gmres,
including a result by Vital \cite{Vital1990}, which follows directly from the containment of the single-vector
Krylov subspace in the block Krylov subspace,
\begin{align*}
	\max_{i=1,\ldots , p}\,\min_{\bt\in\K_m(\bA,\bR_0)}\norm[auto]{\bb^{(i)} - \bA\paren[auto]{(}{)}{\bx_0^{(i)}+\bt}} 
	\leq 
	\max_{i=1,\ldots , p}\,\min_{\bt\in\cK_m(\bA,\br_0^{(i)})}\norm[auto]{\bb^{(i)} - \bA(\bx_0^{(i)} + \bt)}.
\end{align*}

The same subspace containment can be used to show.that for all $i=1,2,\cdots, p$,
\begin{align*}
	\min_{\bt\in\K_m(\bA,\bR_0)}\norm[auto]{\bb^{(i)} - \bA(\bx_0^{(i)}+\bt)} 
	\leq 
	\min_{\bt\in\cK_m(\bA,\br_0^{(i)})}\norm[auto]{\bb^{(i)} - \bA(\bx_0^{(i)} + \bt)}.
\end{align*}

The subspaces underlying \rgmres and block \rgmres satisfy the same containment relationships.
Thus we have
\small
\begin{align*}
	\max_{i=1,2,\ldots p}\min_{\bs\in\cU\atop \bt\in\K_m\paren[auto]{(}{)}{\prn{\bI-\bPhi}\bA,\widehat{\bR}}}& \norm[auto]{\bb^{(i)}-\bA(\widehat{\bx}_0^{(i)} + \bt + \bs)}  
	\\
	\leq& 
  	\max_{j=1,2,\ldots, p}\,\min_{\bs\in\cU\atop \bt\in\cK_m\paren[auto]{(}{)}{\prn{\bI-\bPhi}\bA,\widehat{\br}^{(i)}}} \norm[auto]{\bb^{(i)}-\bA(\widehat{\bx}_0^{(i)} + \bt + \bs)}
\end{align*}
\normalsize
Furthermore, we have that
\begin{align*}
	\min_{\bs\in\cU\atop \bt\in\K_m\paren[auto]{(}{)}{\prn{\bI-\bPhi}\bA,\widehat{\bR}}}\norm[auto]{\bb^{(i)}-\bA(\widehat{\bx}_0 + \bt + \bs)}
	\leq
	\min_{\bs\in\cU\atop \bt\in\cK_m\paren[auto]{(}{)}{\prn{\bI-\bPhi}\bA,\widehat{\br}^{(i)}}}\norm[auto]{\bb^{(i)}-\bA(\widehat{\bx}_0 + \bt + \bs)}.
\end{align*}
In addition, it should be noted that the polynomial approximation interpretation of \gmres has been extended to the block case,
whereby it has been observed that this can be generalized to \emph{matrix-valued polynomials} in the block case
\cite{Simoncini.Conv-Block-GMRES}.  

It should be noted that any per iteration gains realized by using a block method need to be weighed against the additional cost.
Each iteration of a block method requires more FLOPS than the non-block variant, but this comes with the possibility of
accelerated convergence.  Previous researchers have demonstrated that the additional expense of moving to a block method (as
measured in data movement metrics) is only marginally greater than that of its single-vector counterpart. 
We explore this
advantage in \Cref{section.num}. See also the recent results on block methods
on GPUs \cite{Boman.Higgins.Szyld.23}.

\section{Numerical results}\label{section.num}

We have described the original implementation of block \gcrodr, with versions of in \matlab \cite{GCRODR-matlab-code.2016} and a
fully deployed implementation in the \belos package of Sandia's \trilinos Project \cite{PS.GCRODR-Code.2011}.

One point which must be discussed is how to compare the performance of a block \gcrodr to algorithms that execute only a matrix--vector product per
iteration. Block methods have a different dominant core operation in the iteration, the block $p$ matvec.
However, the block $p$ matvec does not cost $p$ times as much as a single standard matvec.
Thus we present two sets of experiments.  One set, shown in \Cref{section.small-experiments}, are all performed in \matlab to demonstrate
characteristics of algorithm performance for small-scale problems.  The second set of experiments,
shown in \Cref{section.data-movement}
are performed
in \trilinos, and demonstrate performance characteristics of the core operations of block \gcrodr
for very large, sparse matrices.
After each cycle, harmonic Ritz vectors are used to build the recycled subspace.

\subsection{Small-scale convergence experiments}\label{section.small-experiments}
The experiments in this section were performed on a Macbook Pro with a 3.1 GHz Dual-Core Intel Core i5 processor and 8 GB of 2133
MHz DDR3 main memory. We demonstrate timing comparisons for performing single and block $p$ matvecs for computing the action of a
sparse matrix $\bA$ on equal numbers of vectors.

\begin{figure}[htb]
\includegraphics[scale=0.55]{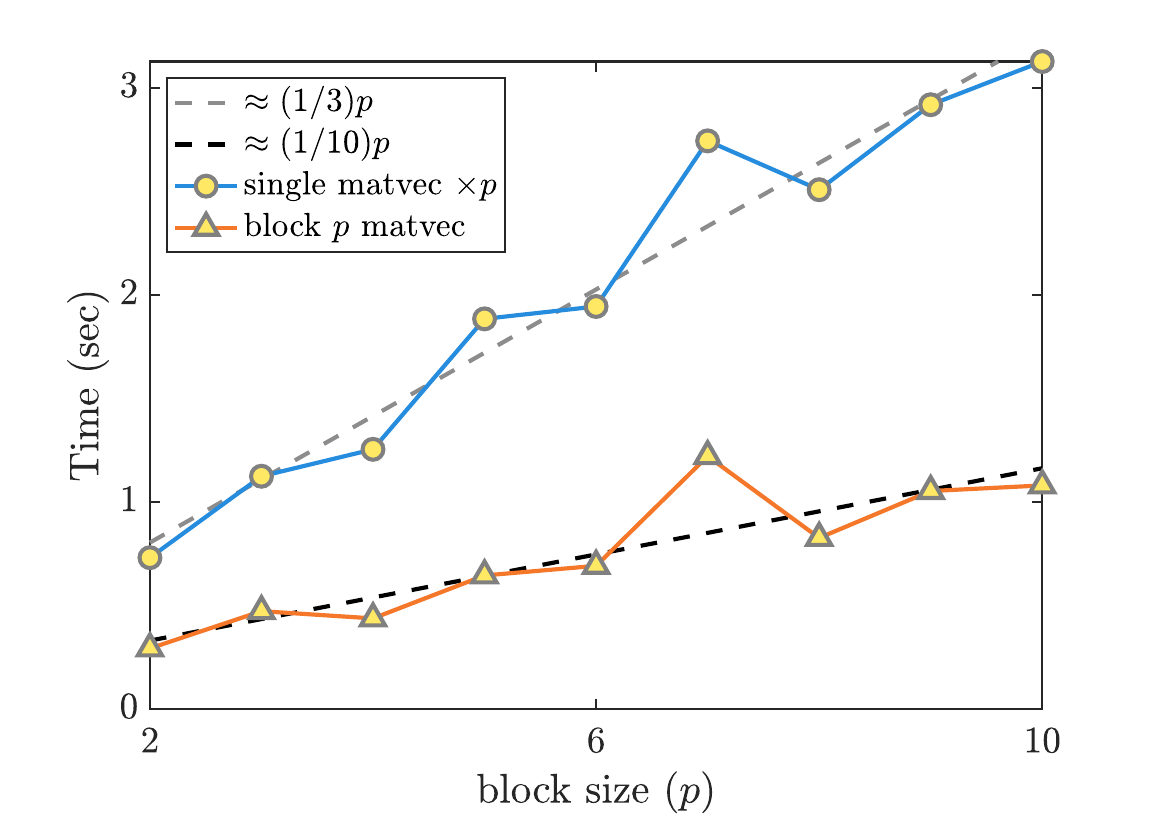}
\caption{Comparison of the time taken to perform $p$ single matrix-vector products for the sparse matrix \texttt{sherman5} of order 3312
	from \cite{2011},  for the single matrix-vector
product versus the block size $p$ block $p$ matvec for different values of $p$, the block size.  
	The \textbf{dashed lines} show the linear least-squares fits to the data from each experiment. Each experiment was
	performed $10^5$ times.
\label{figure.block-vs-single-same-num}}
\end{figure}

\begin{figure}[htb]
\includegraphics[scale=0.55]{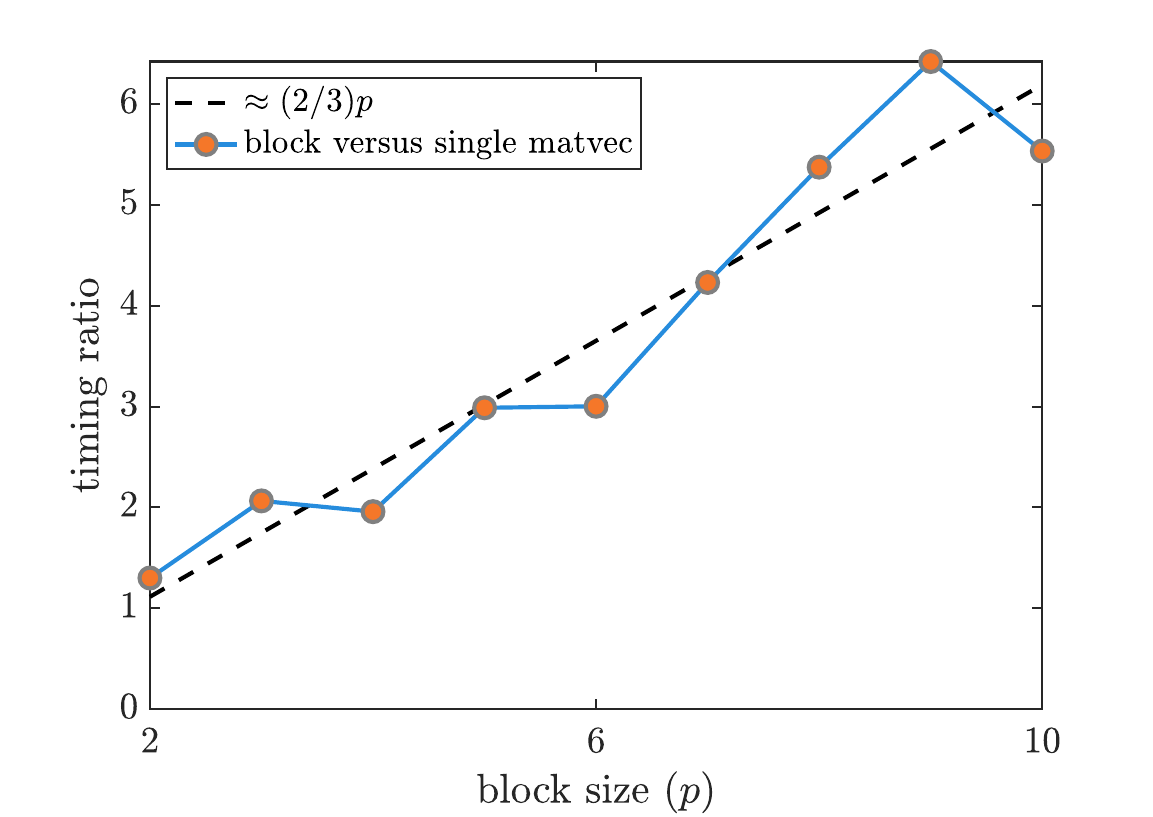}
\caption{
This figure shows the ratio of time taken for the block matrix-matrix products over the time taken for matrix-vector products for
	the \texttt{sherman5} matrix from \cite{2011} for
	different block sizes $p$.  As predicted, the block $p$ matvec is not $p$ times as expensive.  The \textbf{dashed
	line} shows the linear least-squares fit to the experiment data, showing that a block $p$ matvec
	is roughly $\frac{2}{3}p$ times as expensive as a single matrix-vector product, for this matrix.  Each experiment was
	performed $10^5$ times.}
\label{figure.block-vs-single-same-iter}
\end{figure}

In \Cref{figure.block-vs-single-same-num}, we see that the block $p$ matvec is able to outperform the single matrix-vector product, in \matlab.  However, we are more
interested in an per-iteration performance comparison of a block Krylov subspace method versus a single-vector Krylov subspace method.  We pose the
question, do the benefits of convergence in fewer iterations outweigh the increased number of floating-point
operations of the
block block $p$ matvec?  In \Cref{figure.block-vs-single-same-iter}, we compare the time taken to perform matrix-vector products with the time taken to compute block $p$ matvecs.
We see that, though block $p$ matvecs are more expensive to compute than the single-vector variety, they are not $p$ times as
expensive.

\begin{figure}[htb]
\includegraphics[scale=0.5]{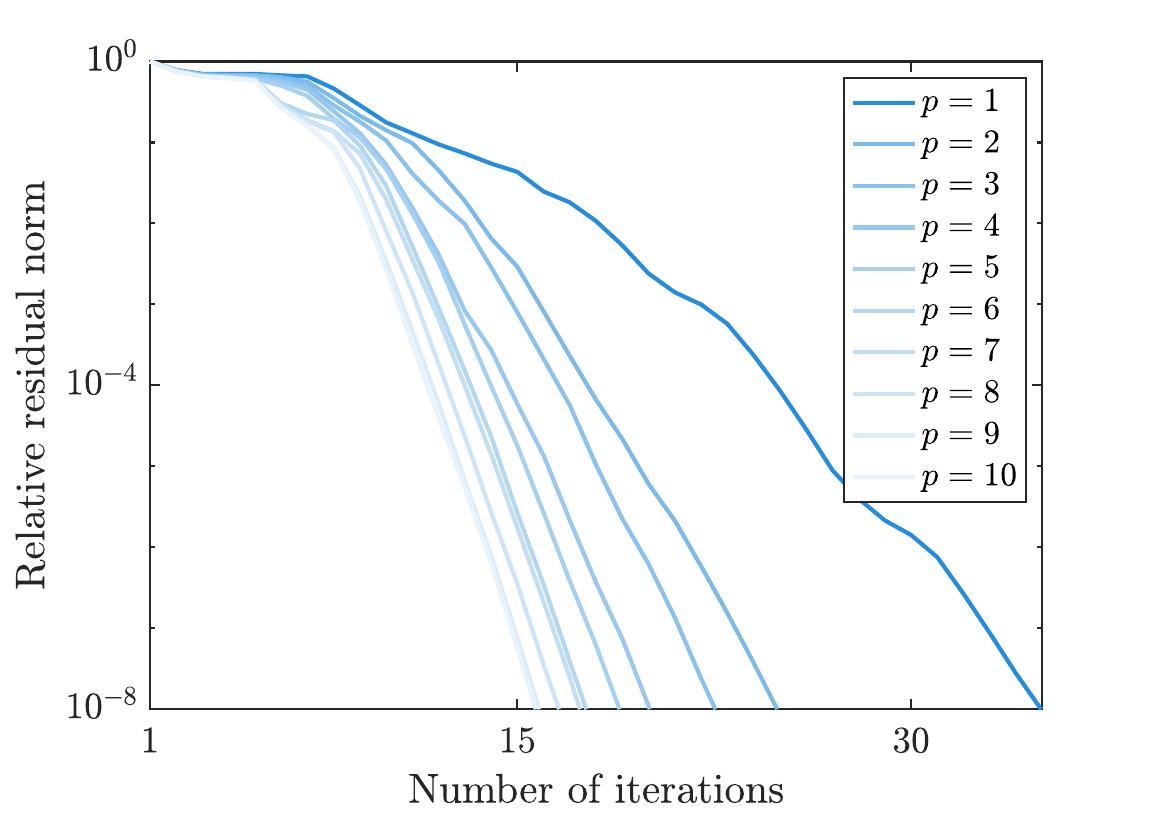}
	\caption{Performance of block recycled \gmres ($m=100$; $k=50$) for the \texttt{sherman5} matrix from \cite{2011} preconditioned with \texttt{ILU(0)} as we increase the
number of right-hand sides.  The first right-hand side is packaged with the matrix; the others are generated using \texttt{rand()}.
Convergence is measured by computing the residual norm associated to the first right-hand side.}
\label{fig.bgcrodr-multiplerhs-sherman5}
\end{figure}

In \Cref{fig.bgcrodr-multiplerhs-sherman5}, we test the code's convergence properties as we increase the number of right-hand
sides.  As is predicted by the underlying theory for this problem, the increased number of right-hand sides generates a richer
space from which to select our approximation updates and from which to recycle, though the marginal benefit decreases for each
additional right-hand side.

We extracted matrices from seven consecutive iterations of a \tramonto Newton iteration from the \texttt{POLY\_CMS\_1D} test problem \cite{dft-poly1-cms-1d}.  For each iteration, we precondition using \texttt{ILU(0)}.
We compare the performance of \gmres, block \gmres with 3 right-hand sides, \rgmres, and our block \rgmres algorithm on all 7 systems.  In the case
of the block methods, one right-hand side generated by the \tramonto package, and the two additional right-hand sides were random, generated using \matlab's \texttt{rand()}.  In the case of these Newton iterations for these relatively small systems
(dimension~$\approx~14000$), convergence for the preconditioned system is fast enough that we are able to recycle the entire Krylov subspace when running
\rgmres and block \rgmres for the first few systems in the sequence before our total subspace dimension exceeds the chosen
dimension $k$ and we must down-select by computing harmonic Ritz vectors.
\begin{figure}[htb]
\includegraphics[scale=0.32]{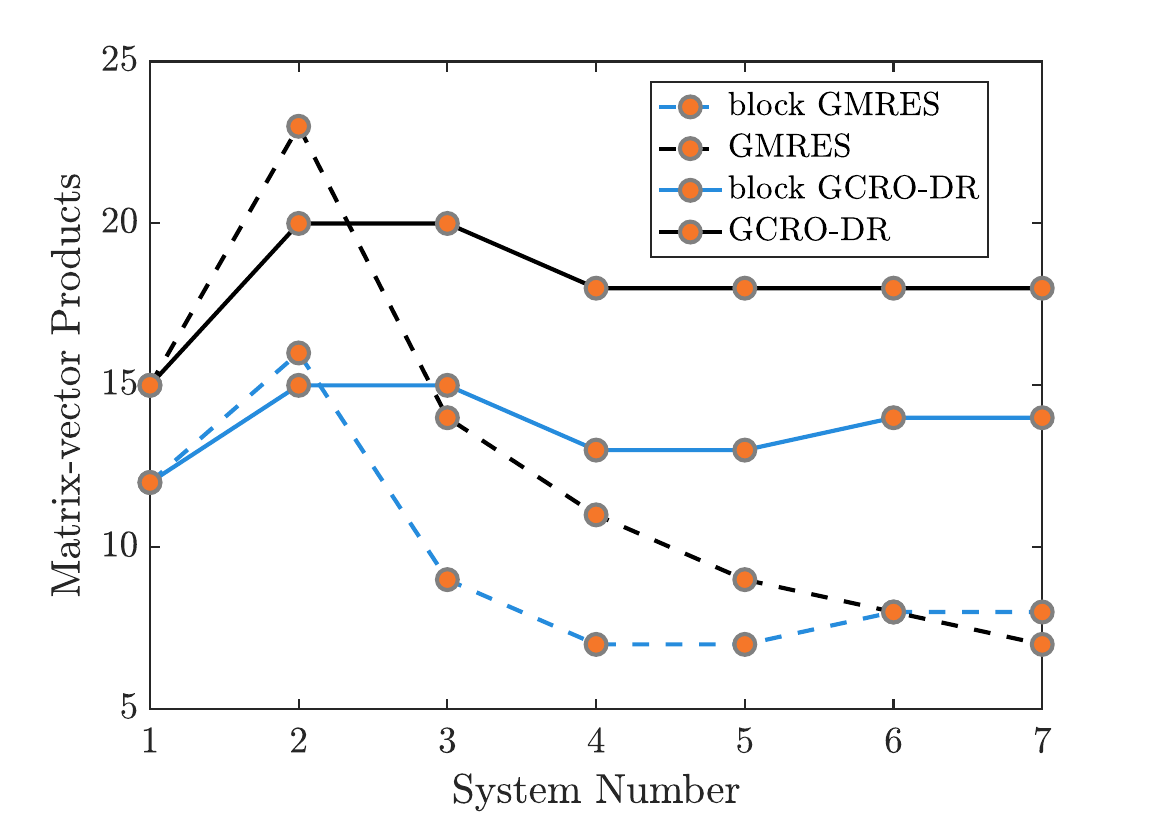}\quad\quad\includegraphics[scale=0.32]{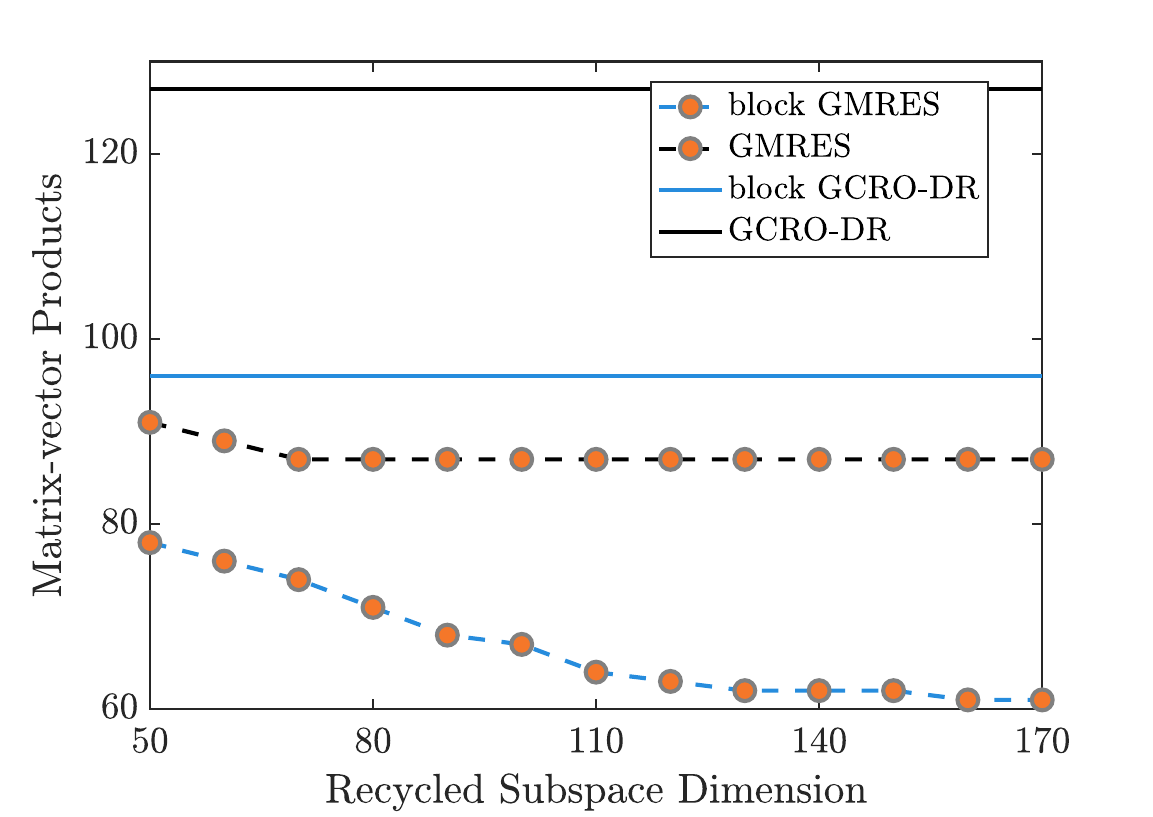}
	\caption{\textbf{Left:} block $p$ matvec count for block recycled \gmres solving linear systems involving the Jacobian for
	7 sequential Newton Steps of fluid DFT problem from \tramonto software package.  ILU(0) preconditioning was used.
	\textbf{Right:} total block $p$ matvec count for various recycled space dimensions   Effectively, as we move to the right
	in this figure, we reach a point at which no recycled subspace information is ever discarded in the experiment.}
\label{fig.precond-bgcrodr-tramonto-7newton}
\end{figure}
In \Cref{fig.precond-bgcrodr-tramonto-7newton}, we plot the number of block $p$ matvecs needed to solve each system.
Observe that for both algorithms, recycling greatly reduces the number of block $p$ matvecs needed to solve later systems.  Furthermore, we
get a per-system reduction when moving from \rgmres to block \rgmres, particularly for systems appearing early in the sequence.
In Figure~\ref{fig.precond-bgcrodr-tramonto-7newton}, we see that for later systems, \gcrodr is able to catch up to the block
method in terms of number of iterations.  This suggests that for some problems for which we use block methods,
the additional expense of recycling may bring the most benefit for the earlier systems.
This can yield a high-quality recycled subspace, and we may then be able to apply single-vector \rgmres
for the rest of the systems.   We also see that, for large enough recycled subspace, we achieve a 30\% reduction in overall
matvecs when moving from single right-hand side \rgmres to block \rgmres with two random right-hand sides.

\subsection{Larger parameter study}\label{section.parameter-study}
In this section, we perform a parameter study of block \gcrodr.  The matrices
we use arise from finite element discretization of the steady-state
convection-diffusion problem, posed on the square $\Omega = [-1,1]^2$,
generated in IFISS \cite{ers07,ers14,ifiss} from the standard
\emph{double-glazing} problem \cite[Example 3.1.4]{ESW-Book.2014} using the
built-in script \texttt{square\_cd} with specified parameter of $10$. The
diffusion coefficient was generated randomly to simulate uncertainty in that
coefficient, from a distribution with mean $0.01$.  This yields matrices of
size $1046529\times 1046529$.  Incomplete-\lu factorization with a drop 
tolerance of $10^{-1}$ was used to generate the preconditioner for each system.
\begin{figure}[htb]
	\centering
	\includegraphics[scale=0.25]{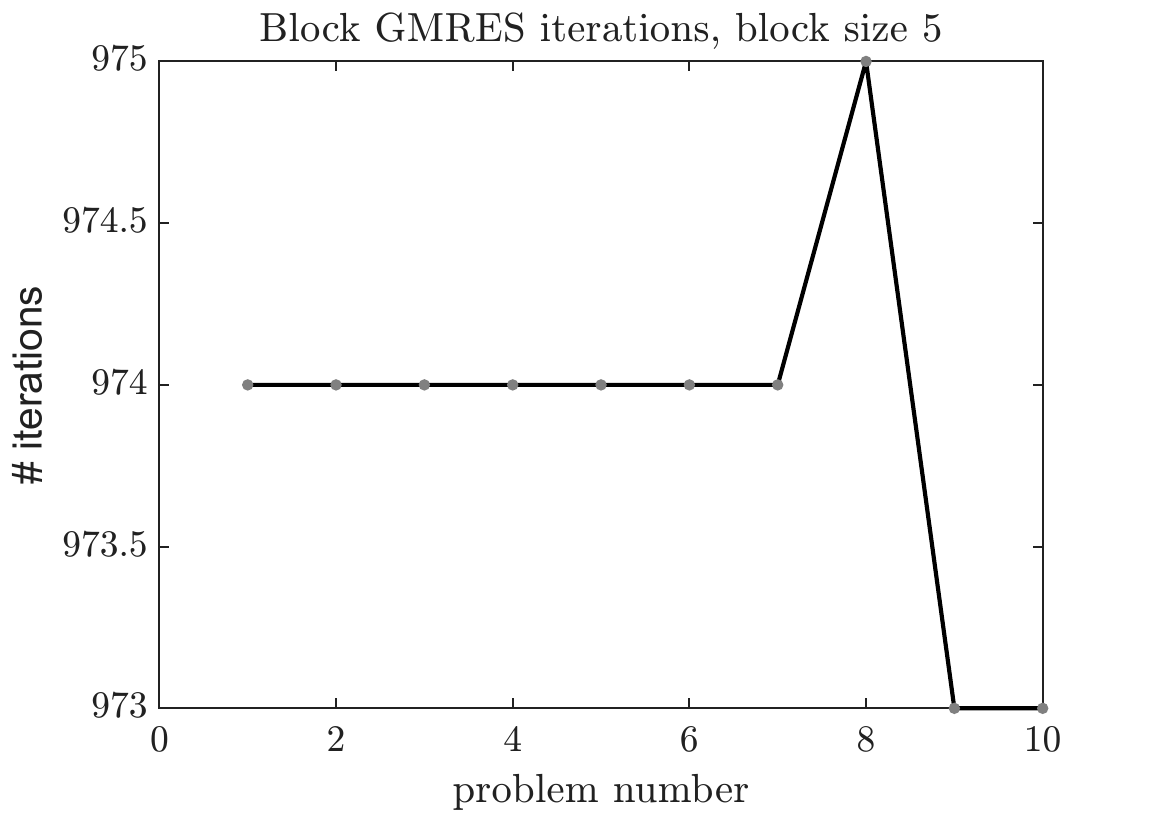}
	\includegraphics[scale=0.25]{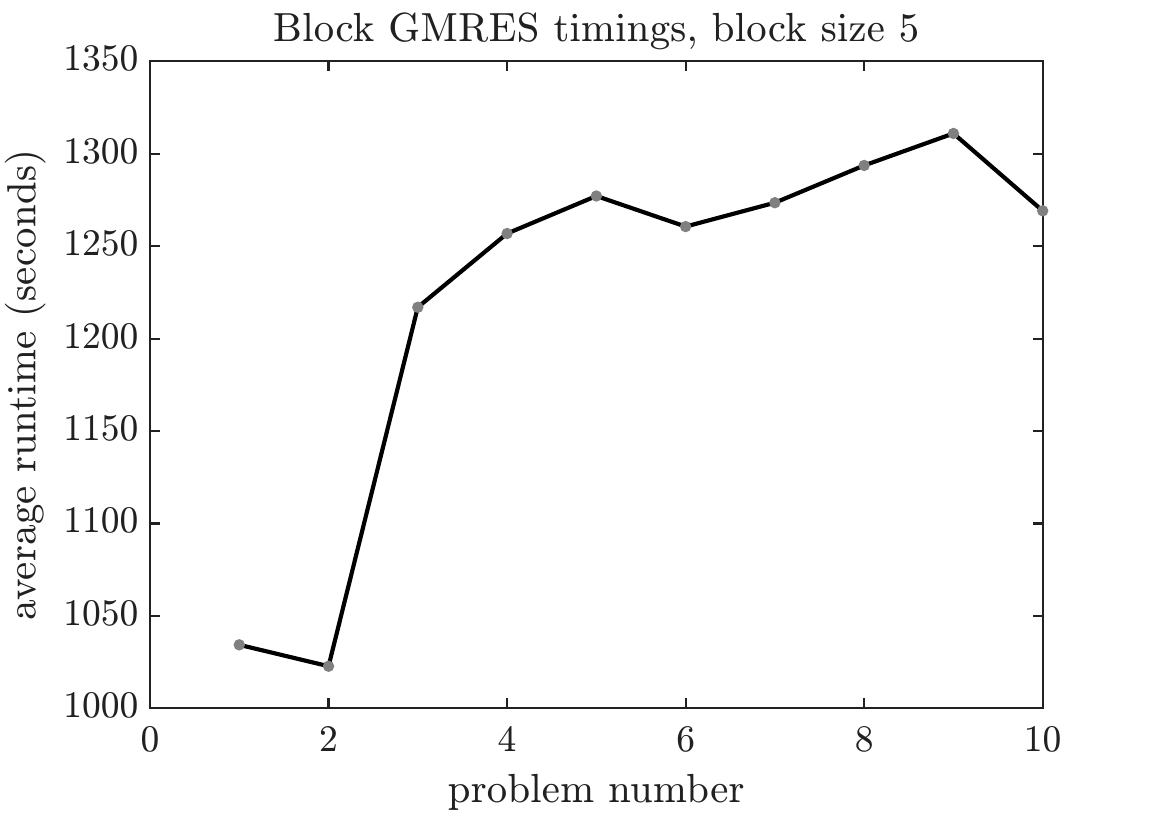}
	\caption{
		\label{figure:block-gmres-cd-tests}
		Iteration count and runtimes for \blgmres for the ten 
		problems described in \Cref{section.parameter-study} for block-size $5$.
	}
\end{figure}

\begin{figure}[htb]
	\centering
	\includegraphics[scale=0.25]{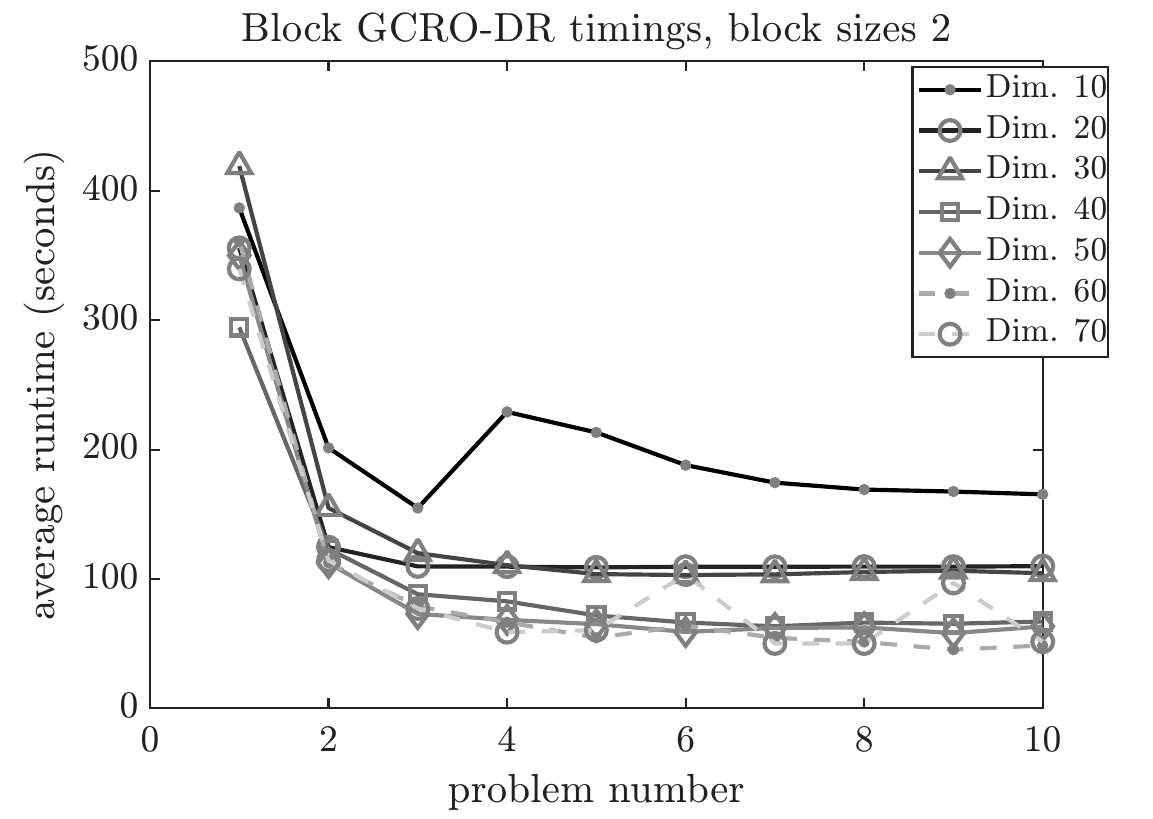}
	\includegraphics[scale=0.25]{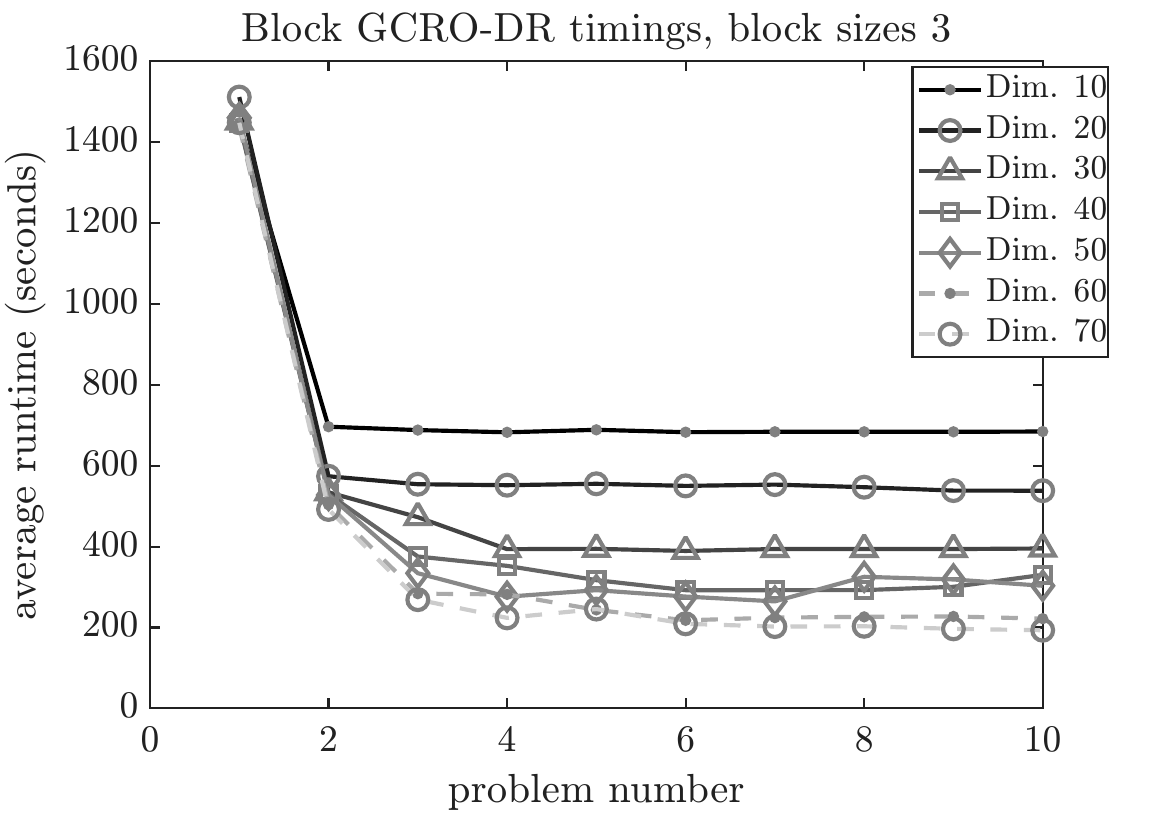}
	\\
	\includegraphics[scale=0.25]{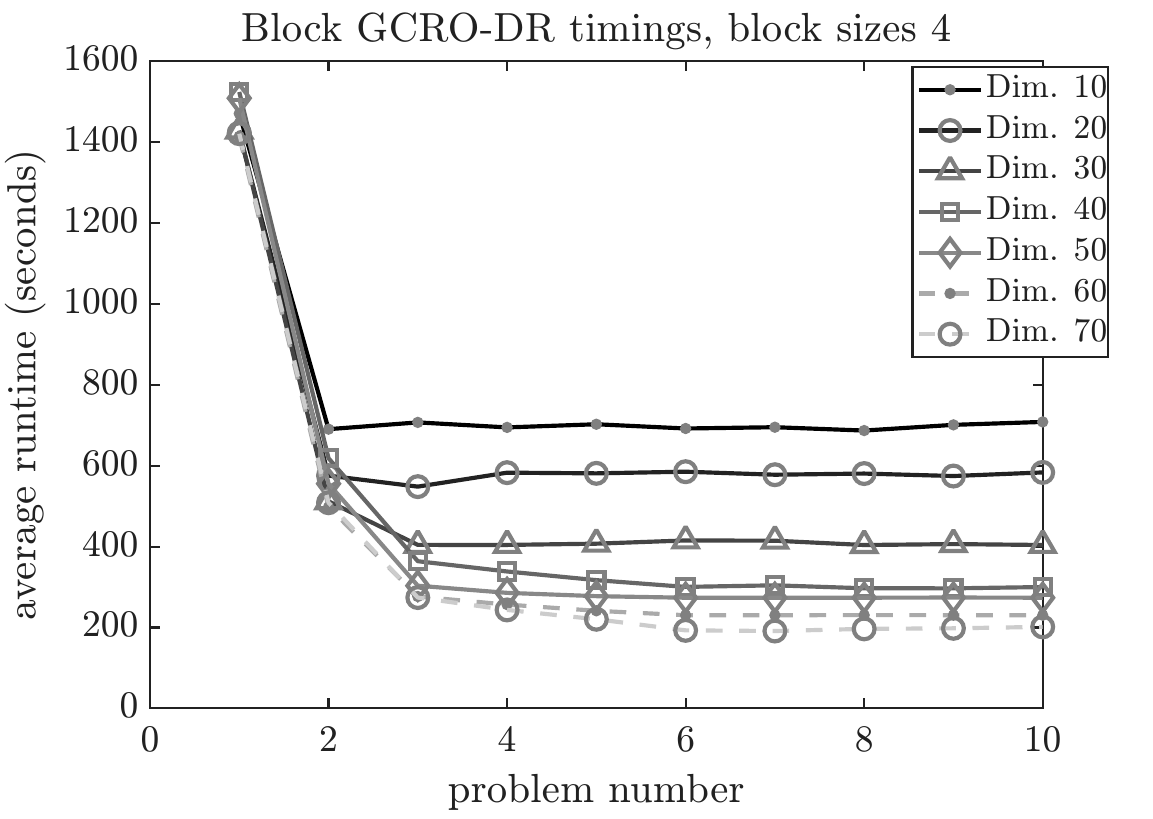}
	\includegraphics[scale=0.25]{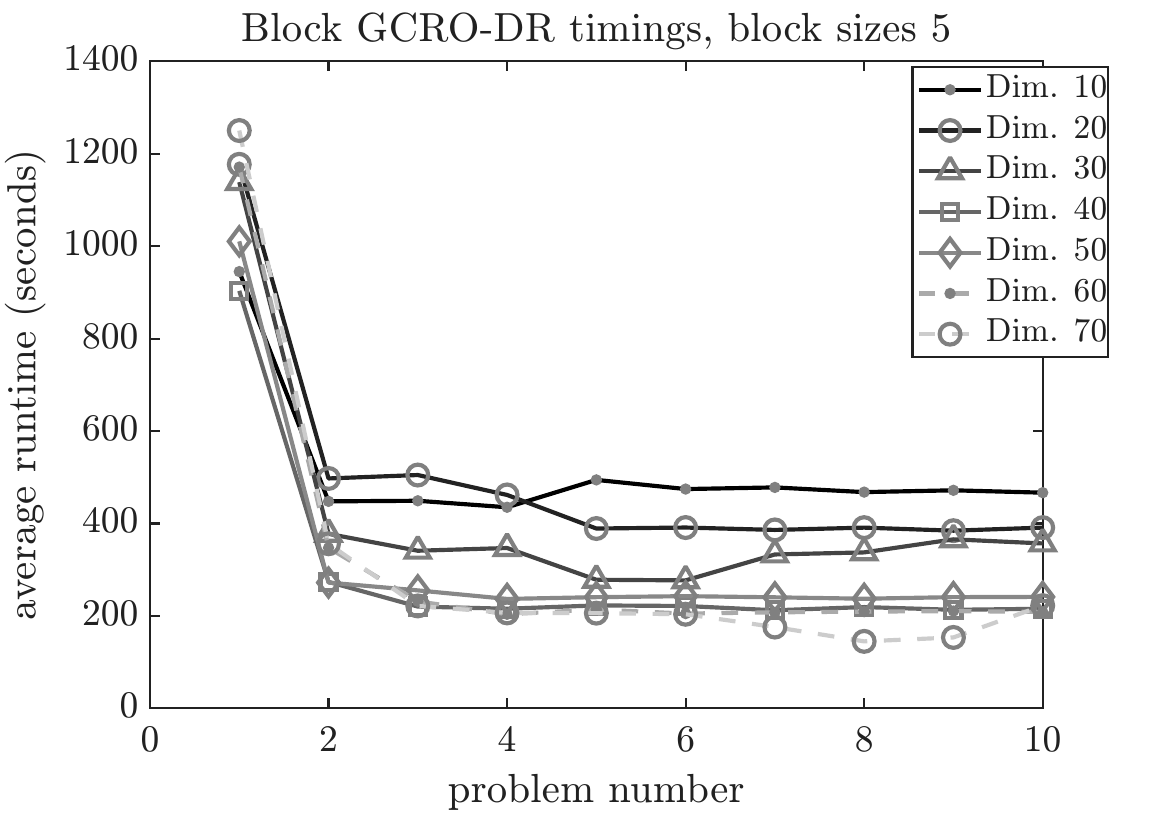}
	\caption{
		\label{figure:block-gcrodr-timings}
		Runtimes for block \gcrodr for the ten 
		problems described in \Cref{section.parameter-study} for block-sizes $2-5$
		and various recycled subspace dimensions.  Cycle-length as adjusted so that 
		storage requirements were roughly constant across different runs.
	}
\end{figure}

Five right-hand sides were generated for each matrix.  We ran our block \gcrodr
\matlab code on this sequence for block sizes $2, 3, 4, 5$ and with recycled
\begin{figure}[htb]
	\centering
	\includegraphics[scale=0.25]{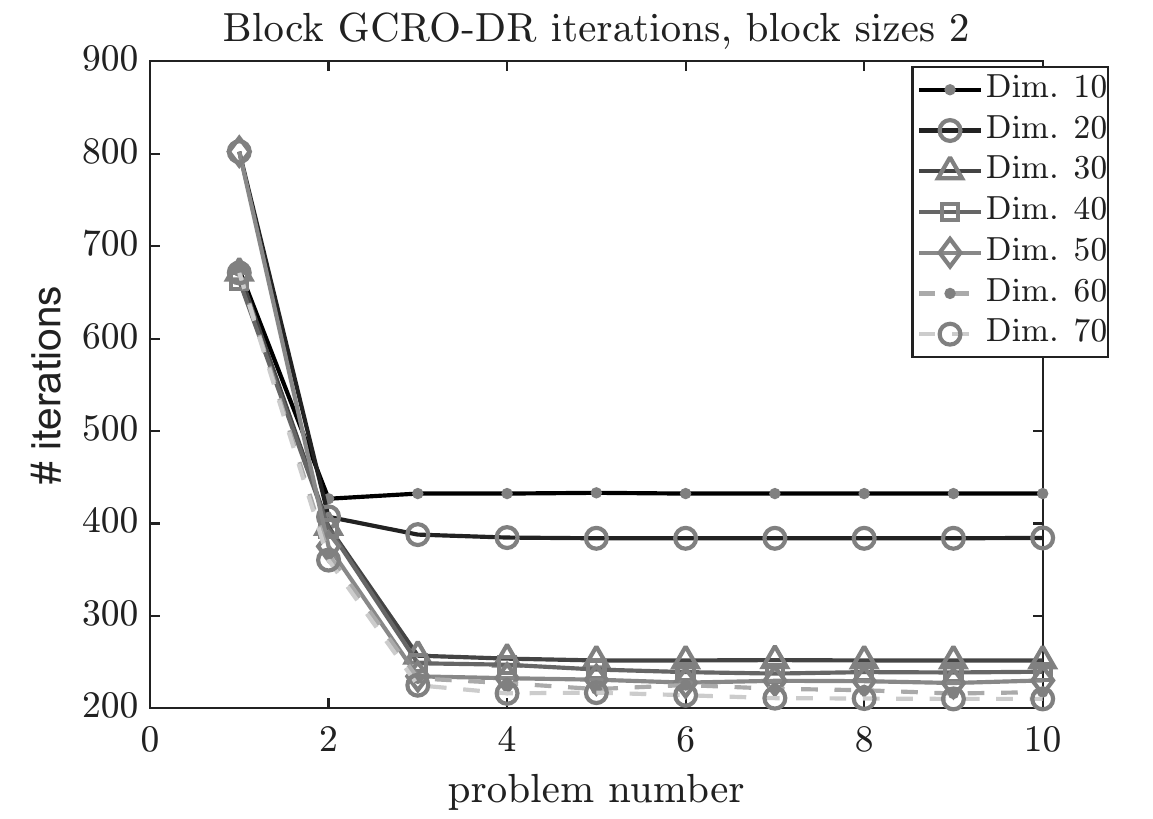}
	\includegraphics[scale=0.25]{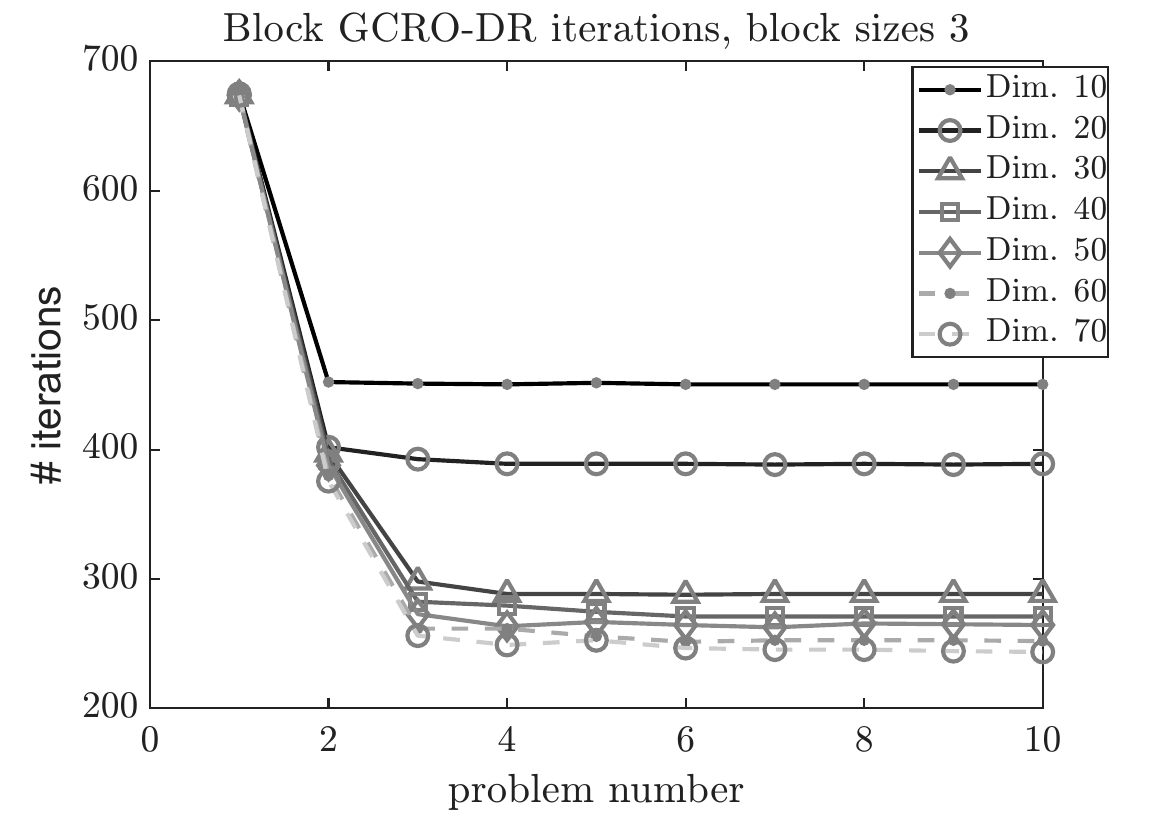}
	\\
	\includegraphics[scale=0.25]{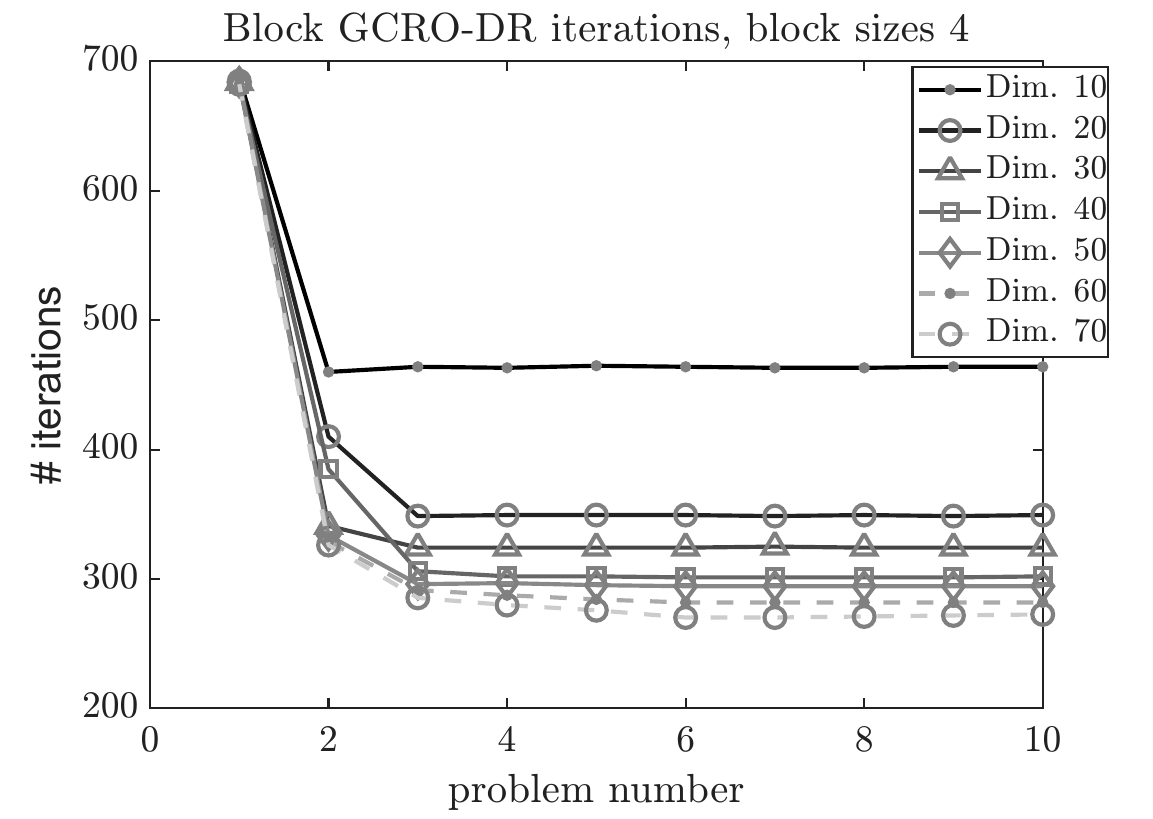}
	\includegraphics[scale=0.25]{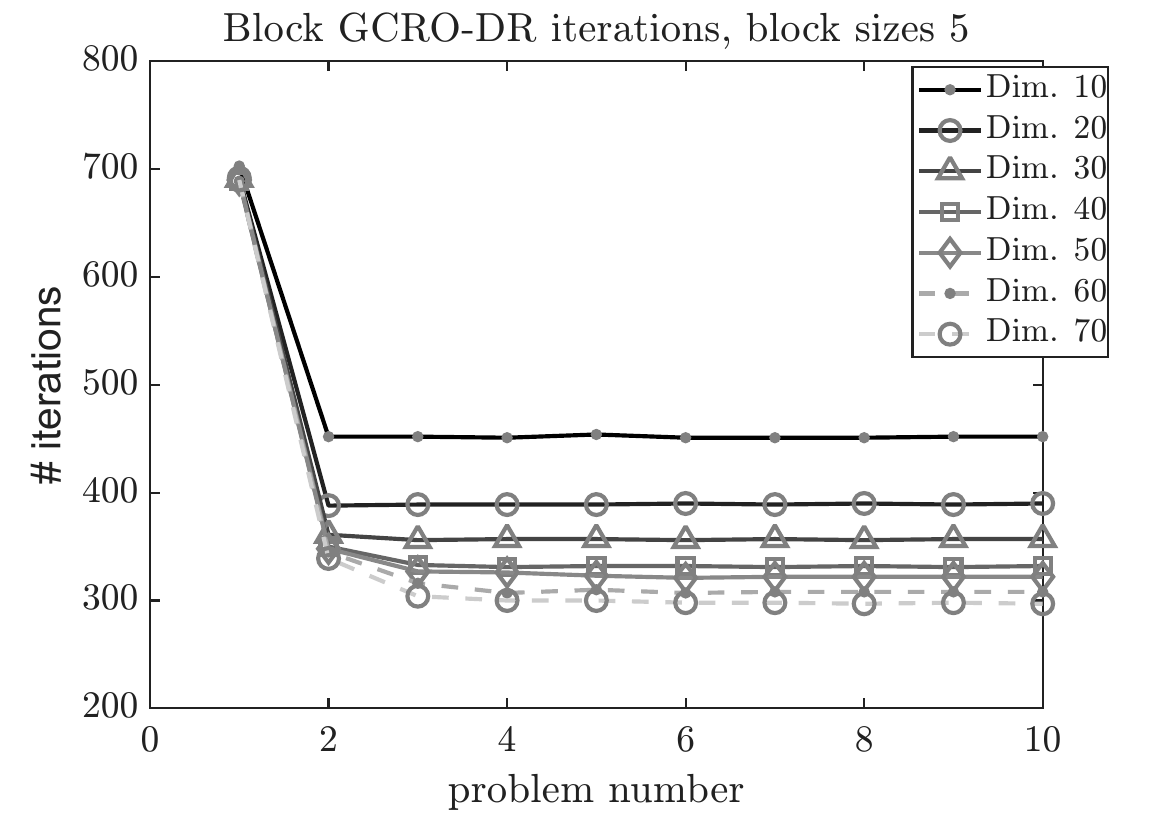}
	\caption{
		\label{figure:block-gcrodr-iter}
		Iteration counts for block \gcrodr for the ten 
		problems described in \Cref{section.parameter-study} for block-sizes $2-5$
		and various recycled subspace dimensions.  Cycle-length was adjusted so that 
		solution space dimensions were roughly constant across different runs.
	}
\end{figure}
subspace dimensions $10, 20, 30, 40, 50, 60, 70$. For each recycled subspace
dimension, the cycle length was adjusted so that the augmented Krylov subspace
dimension was as close to $200\times \mathtt{blockSize}$ as possible, to yield a fairer 
comparison in performance for roughly equidimensional solution subspaces.
For this particular problem, we found that recycling was more effective when
selecting the largest harmonic Ritz values at the end of each cycle, rather
than the smallest. These tests were run on \macbook Pro with an Apple M4 Max
processor and 36GB of shared memory.  We used \matlab version R2025a. 
Note: 
\emph{
	Due to the size of these problems and their run-times, we did not
	pursue the usual practice of running each test multiple times to get
	average timings. Instead we report the timings from a single run per
	experiment.
}
\Cref{figure:block-gmres-cd-tests} shows iterations and timings 
for \blgmres with no recycling for five right-hand sides. \Cref{figure:block-gcrodr-iter} shows 
iteration counts for block \gcrodr for the sequence of ten problems for block sizes $2-5$.
\Cref{figure:block-gcrodr-timings} shows timings for the same experiments.
We observe that the marginal benefit of adding more vectors to the recycled subspace 
decreases as we allow for larger recycled subspaces.  

These experiments illustrate that, although recycling methods can be used as
block-box solvers with the default settings, it behooves the end-user to
experiment and understand what works best for a specific application problem to
get the best performance at the lowest cost.

\begin{table}[htb]
\caption{Fifteen different matrices downloaded from the University of Florida Sparse Matrix Library \cite{DH.2011} with various sparsity patterns.  The matrix \texttt{qcdRealPart} is the real part of the matrix \texttt{conf5-0-4x4-10} from \cite{DH.2011}}
\label{table.uniFloridaMatInfo}
\begin{tabular}{l|r|r|c}
Name & Dimension & \# Non-zeros & Sparsity \\
\hline
\texttt{Freescale1} & $3428755$ & $17052626$ & $0.0001\%$ \\
\texttt{CoupCons3D} & $416800$ & $17277420$ & $0.0099\%$ \\
\texttt{rajat31} & $4690002$ & $20316253$ & $0.0001\%$ \\
\texttt{FullChip} & $2987012$ & $26621983$ & $0.0003\%$ \\
\texttt{cage14} & $1505785$ & $27130349$ & $0.0012\%$ \\
\texttt{RM07R} & $381689$ & $37464962$ & $0.0257\%$ \\
\texttt{epb3} & $84617$ & $463625$ & $0.0065\%$ \\
\texttt{qcdRealPart} & $49152$ & $1916928$ & $0.0793\%$ \\
\texttt{crashbasis} & $160000$ & $1750416$ & $0.0068\%$ \\
\texttt{Hamrle3} & $1447360$ & $5514242$ & $0.0003\%$ \\
\texttt{HV15R} & $2017169$ & $283073458$ & $0.0070\%$ \\
\texttt{lung2} & $160000$ & $1750416$ & $0.0068\%$ \\
\texttt{ML\_Geer} & $1504002$ & $110686677$ & $0.0049\%$ \\
\texttt{pre2} & $659033$ & $5834044$ & $0.0013\%$ \\
\texttt{twotone} & $120750$ & $1206265$ & $0.0083\%$ 
\end{tabular}
\end{table}
\begin{figure}[htb]
\includegraphics[scale=0.25]{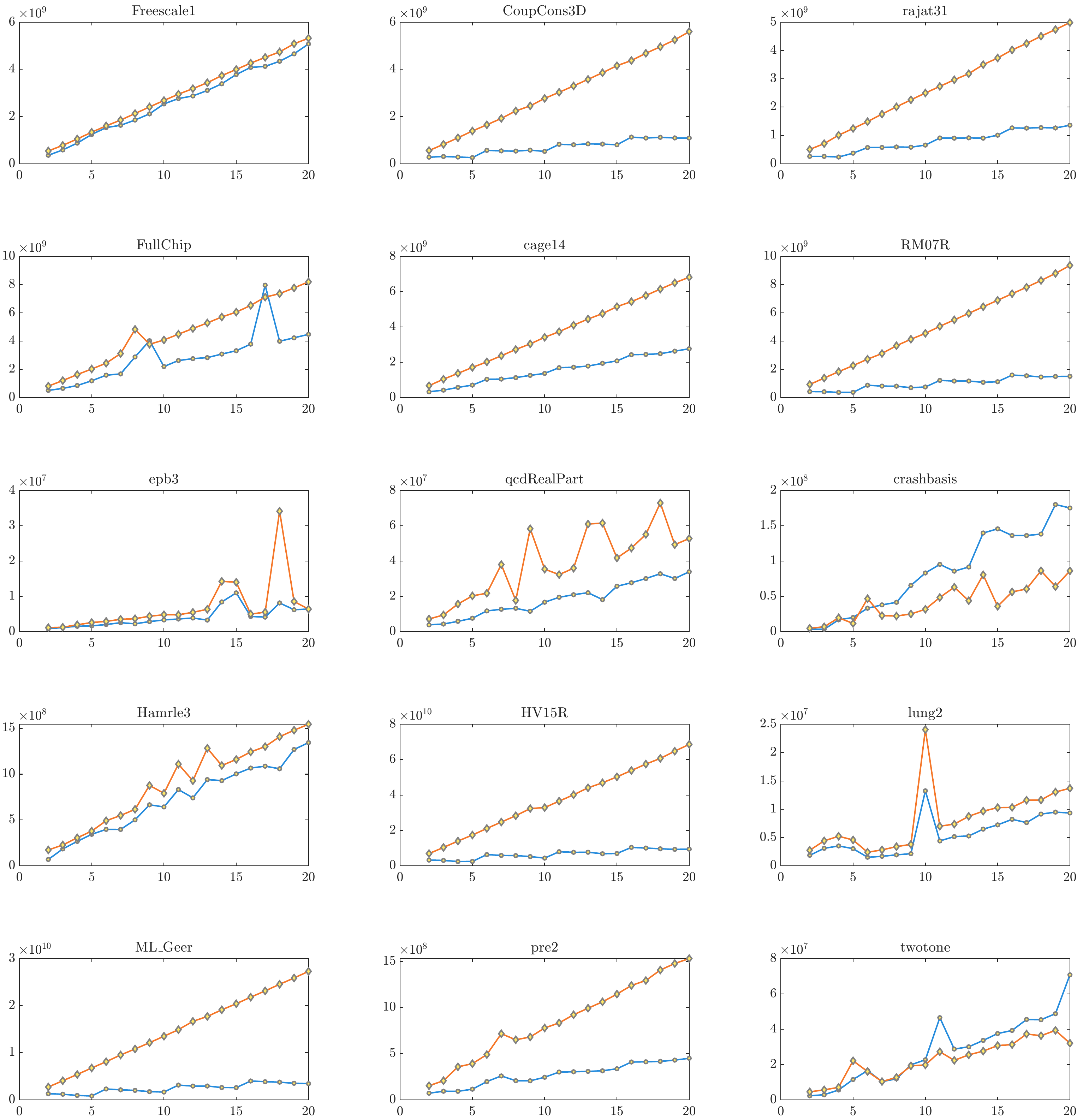}
\caption{\label{figure.cacheMissesFloridaComp} For the fifteen different matrices from \Cref{table.uniFloridaMatInfo}, a comparison of cache misses for the matrix applied to different sizes of vector blocks both all-at-once (\textbf{blue lines with yellow circles}) and one-at-a-time (\textbf{red lines with yellow diamonds}).  The horizontal axes show block size and the vertical axes show number of cache misses per 100 executions of the experiment. }
\end{figure}

\begin{figure}[htb]
\includegraphics[scale=0.27]{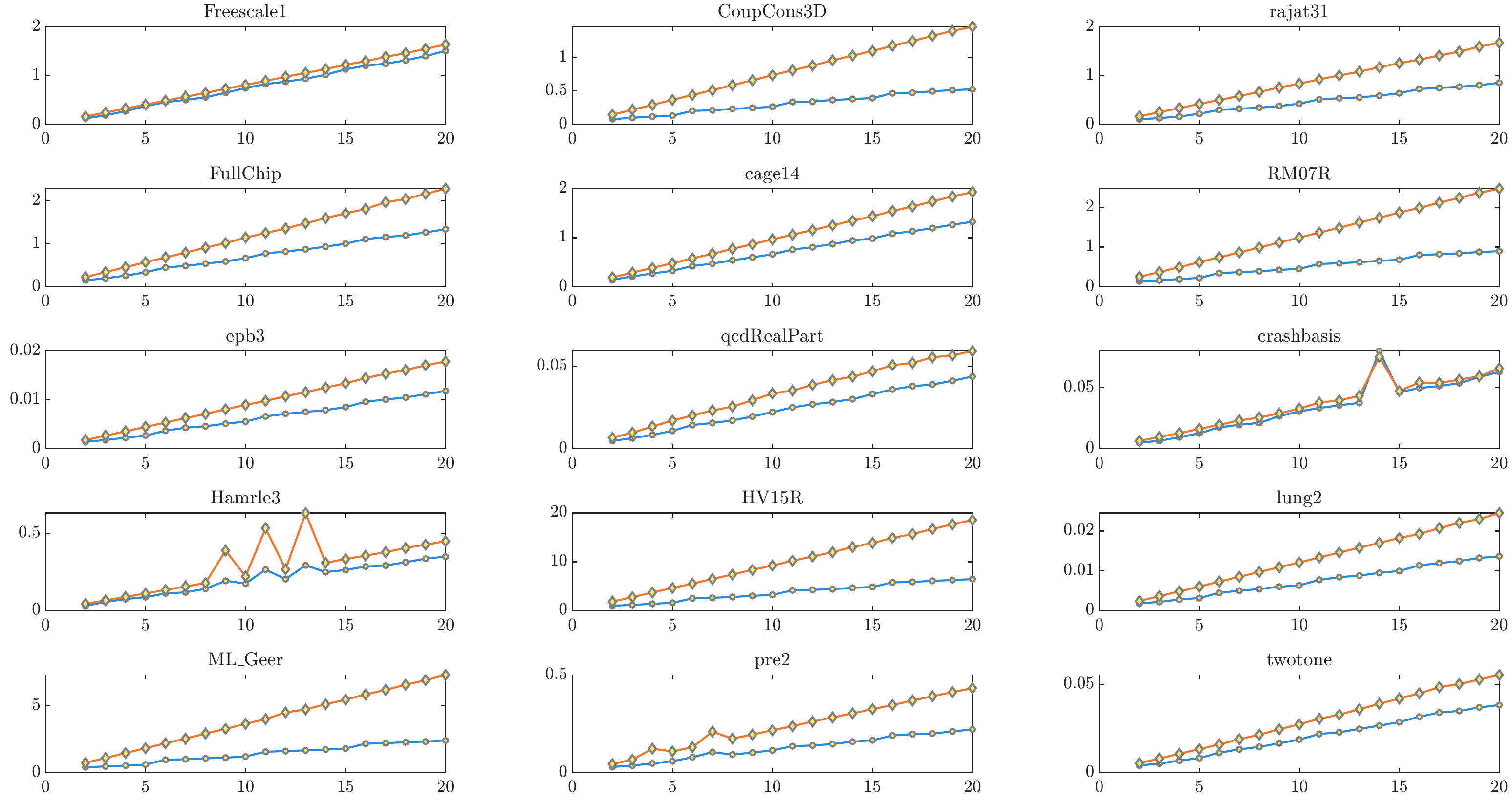}
\caption{\label{figure.avgTimingsFloridaComp} For the fifteen different matrices from \Cref{table.uniFloridaMatInfo}, a comparison timings (in seconds) for the projected matrix applied to different sizes of vector blocks both all-at-once (\textbf{blue lines with yellow circles}) and one-at-a-time (\textbf{red lines with yellow diamonds}).  The horizontal axes show block size and the vertical axes show the time in seconds required for each experiment. Note that in each plot, the vertical axis is
in a different scale.}
\end{figure}

\begin{figure}[htb]
\includegraphics[scale=0.25]{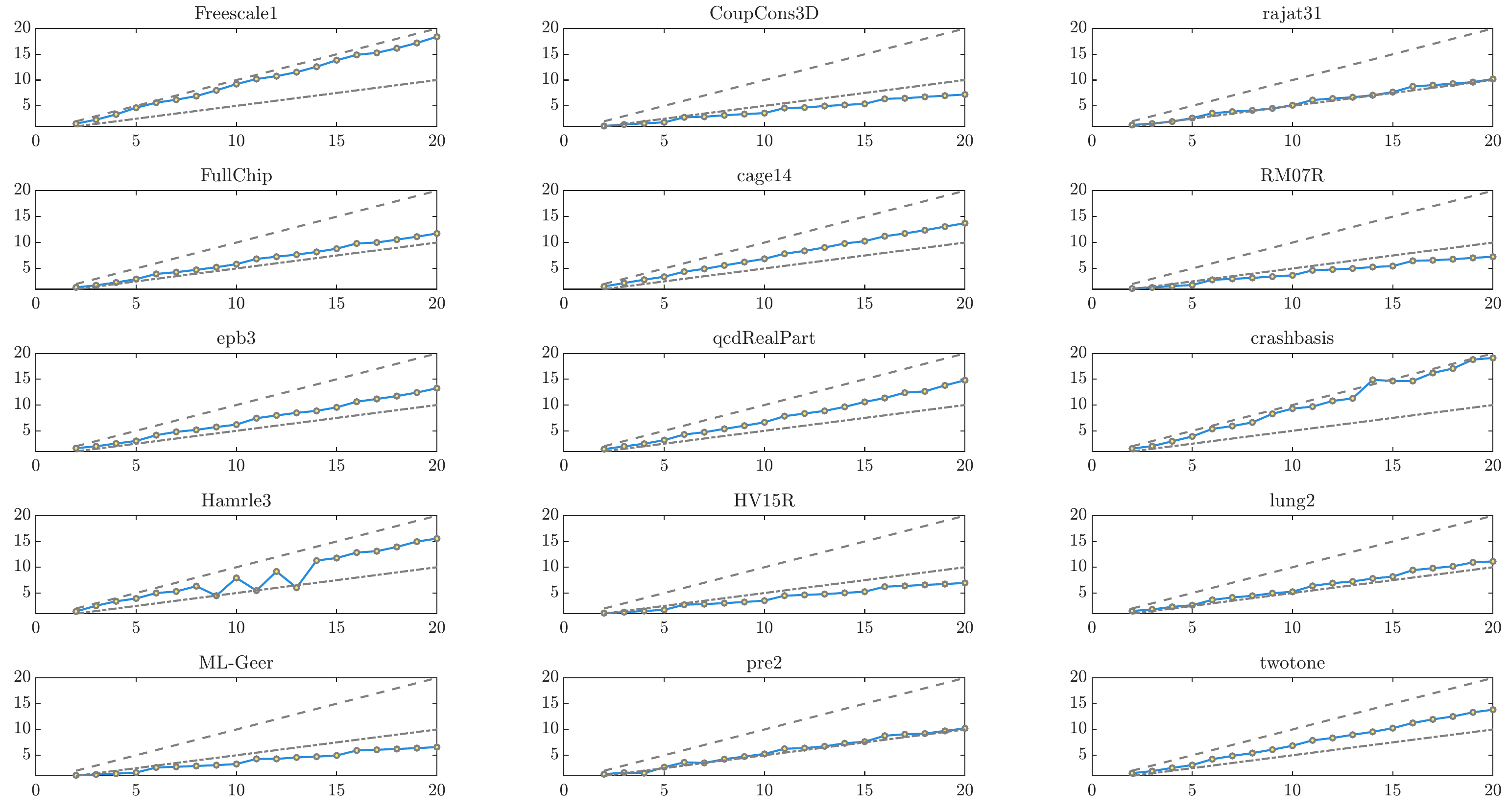}
\caption{\label{figure.blockSingleRatioFloridaComp} For the fifteen different matrices from \cite{DH.2011}, the ratio of the time taken to apply the matrix to different sizes of vector blocks versus to a single vector. The horizontal axes show block size and the vertical axes show the ratio for each block.  For reference, we include references lines for ratio $p$ \textbf{\rm (gray dashed line)} and $\frac{1}{2}p$ \textbf{\rm (gray dot-dashed line)}. }
\end{figure}    

\subsection{Data movement experiments}\label{section.data-movement}

\begin{figure}[htb]
	\includegraphics[scale=0.25]{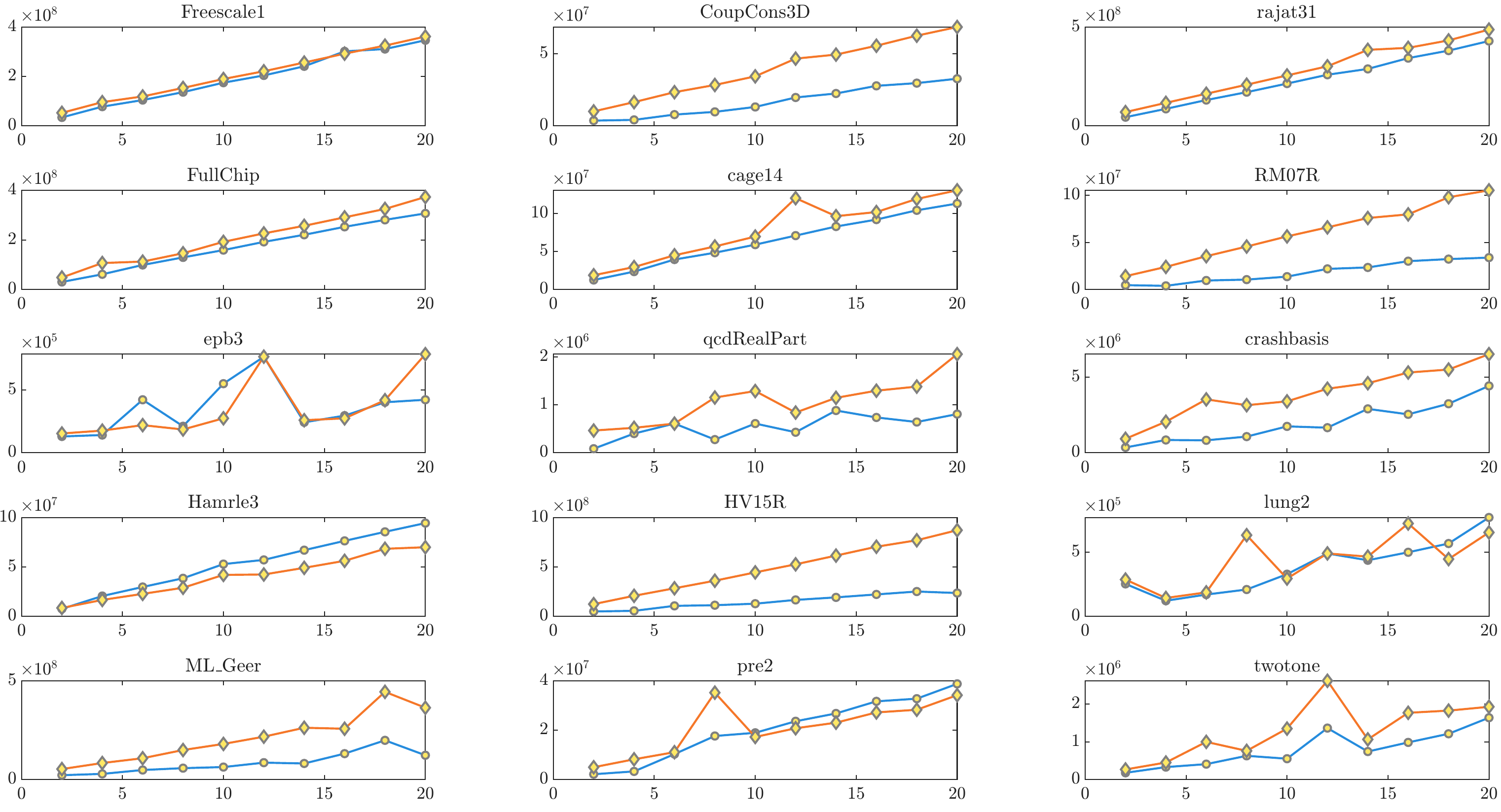}
	\caption{\label{figure.cacheMissesProjectedFloridaComp} For the fifteen different matrices from \Cref{table.uniFloridaMatInfo}, a
	comparison of cache misses for the projected matrix $(\bI - \bPhi)\bA $ applied to different sizes of vector
	blocks both all-at-once (\textbf{blue lines with yellow circles}) and one-at-a-time (\textbf{red lines with yellow diamonds}).  The projector was applied using a
	modified Gram-Schmidt algorithm, and $\bC\in\R^{n\times 100}$.  The horizontal axes show block size and the vertical axes
	show number of cache misses per 100 executions of the experiment. }
\end{figure}

\begin{figure}[htb]
	\includegraphics[scale=0.25]{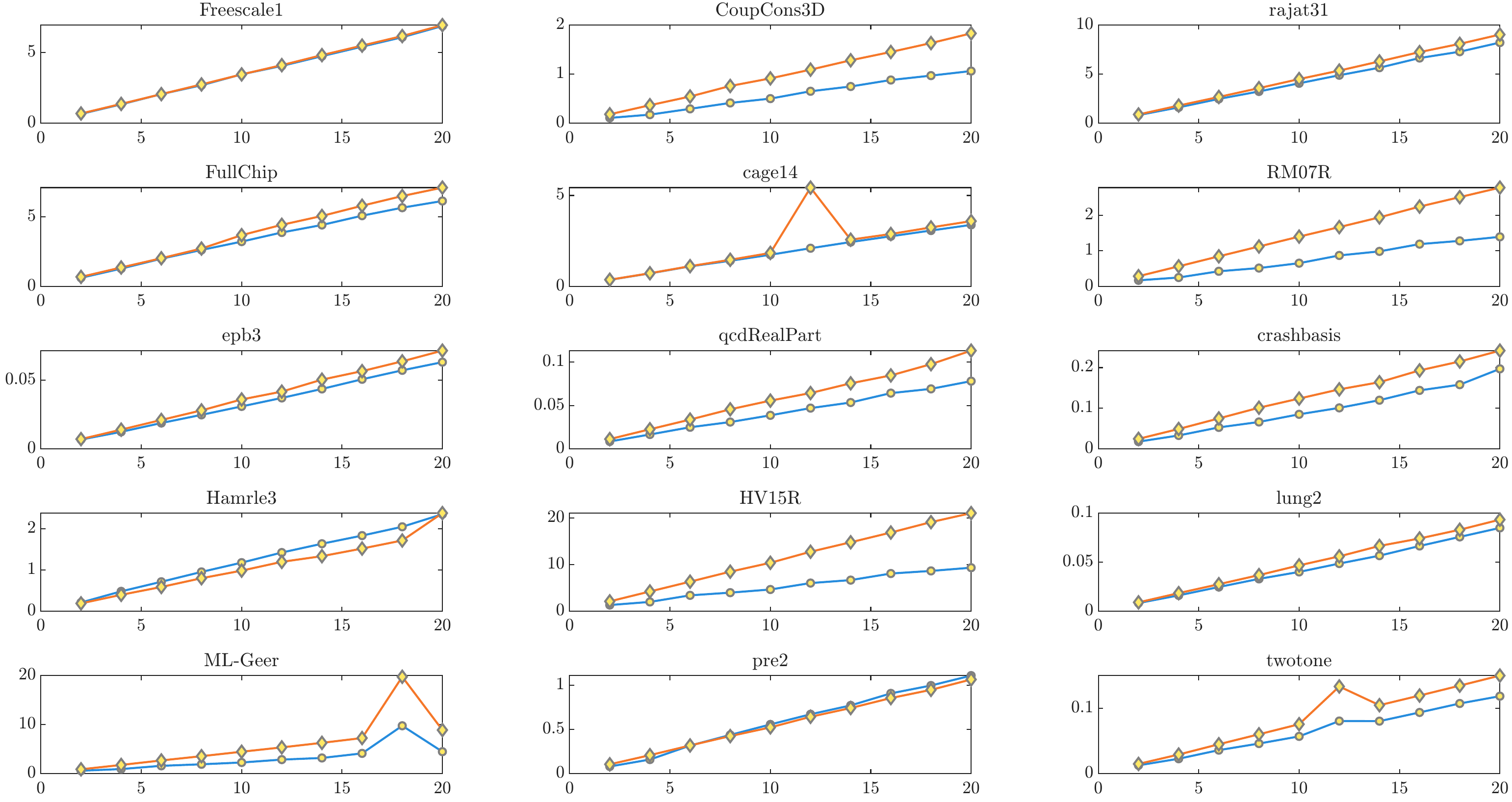}
	\caption{\label{figure.avgTimingsProjectedFloridaComp} For the fifteen different matrices from \cite{DH.2011}, a comparison
	average timings (in seconds) for the projected matrix $(\bI - \bPhi)\bA $ applied to different sizes of vector
	blocks both all-at-once (\textbf{blue lines with yellow circles}) and one-at-a-time (\textbf{red lines with yellow diamonds}).  The projector was applied using a
	modified Gram-Schmidt algorithm, and $\bC\in\R^{n\times 100}$.  The horizontal axes show block size and the vertical axes
	show the average time in seconds required for each experiment. }
\end{figure}

\begin{figure}[htb]
	\includegraphics[scale=0.25]{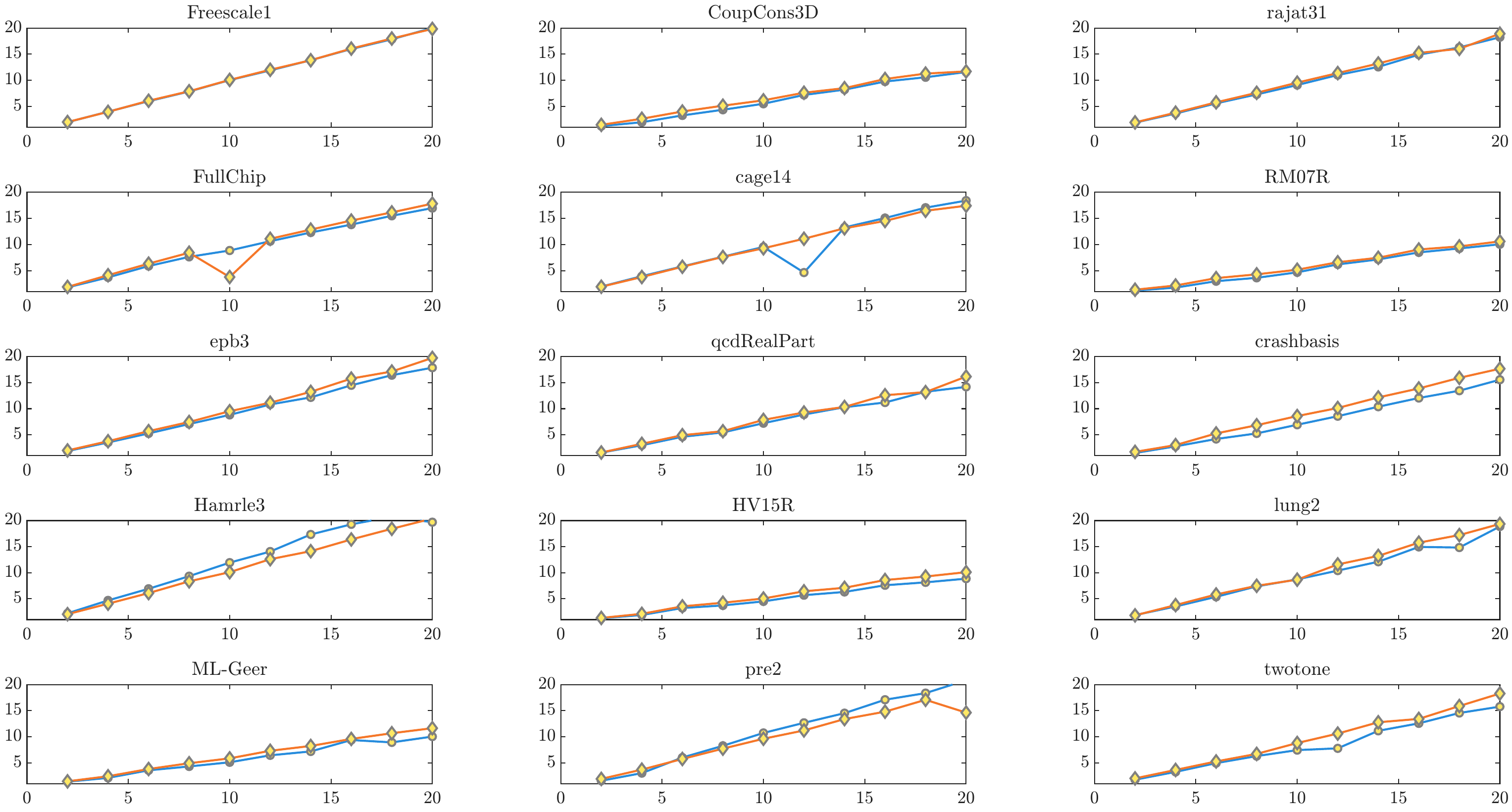}
	\caption{\label{figure.blockProjectedSingleRatioFloridaComp} For the fifteen different matrices from \cite{DH.2011}, the ratio of
	the average time taken to apply the matrix to different sizes of vector blocks versus to a single vector for the projected
	operator $\paren[auto]{(}{)}{\bI - \bPhi}\bA$. The projector is applied using either a modified Gram-Schmidt algorithm
	(\textbf{blue lines with yellow circles}) or multiple passes of classical Gram-Schmidt (\textbf{red lines with yellow diamonds}), and $\bC\in\R^{n\times 100}$.
	The horizontal axes show block size and the vertical axes show the ratio for each block size. }
\end{figure}

In this section, we run a variety of performance tests on matrices of various
sparsity patterns and levels coming both from real applications \cite{DH.2011}
 and from test sets we artificially constructed.
The tests were performed on a shared memory machine with 8 Intel Xeon E7-4870 2.4Ghz processors, each with
10 cores (i.e., a total of 80 CPU cores) and a total of 1 terabyte main memory.  Each processor has 30 megabytes
of L3 cache (3 megabytes per core).  
Each experiment was run as a single core, that is, in serial mode, without threading or MPI. 
Using compiled
\trilinos codes \cite{2011a}, we compare performance of large sparse operators being applied to
blocks of vectors of varying block sizes.  For each matrix $\bA$, we compare
multiplying the matrix times the entire block versus multiplying the matrix times each vector individually. We note that others \cite{doi:10.1137/20M1382064}  have also explored using a block size smaller than the number of right-hand sides for reasons of cache size, although we do not explore that variation here. When multiplying a matrix times a vector individually we store the vector as an Epetra\_Vector object, which consists of a single array of double-precision values that are always contiguous in memory. When multiplying a matrix times a block of vectors, the vectors are stored as an Epetra\_Multivector. An Epetra\_MultiVector may be thought of as a generalization of a dense matrix. This object stores all vectors one after another in a single large array. For each matrix, the experiment
is repeated for the projected operator $\paren[auto]{(}{)}{\bI - \bPhi}\bA$ for subspaces $\cU$
of different dimensions.  Before each test was performed, a block $p$ matvec was executed so that
any prefetching of data into the cache would occur before the start of the test and thus would not interfere with
our measurements.

Performance is measured in multiple ways.  First, each experiment is performed $100$ times and the average
time in seconds for those experiments is taken.  Second, we compiled the Intel Performance Counter Monitor (PCM)
libraries \cite{PCM.2016} and inserted appropriate function calls into our test code,
and these were used to take measurements directly from the processor for
each experiment.  Namely, the PCM allows one to measure bytes read by the processor, the percentage of
cache hits, and the number of cache misses occurring during the experiment.  
In our experiments, we demonstrate that
often, the sparse block $p$ matvec has superior performance when measured in these
cache- and data- related metrics over the sparse matrix-vector product.

\begin{figure}[htb]
\includegraphics[scale=0.25]{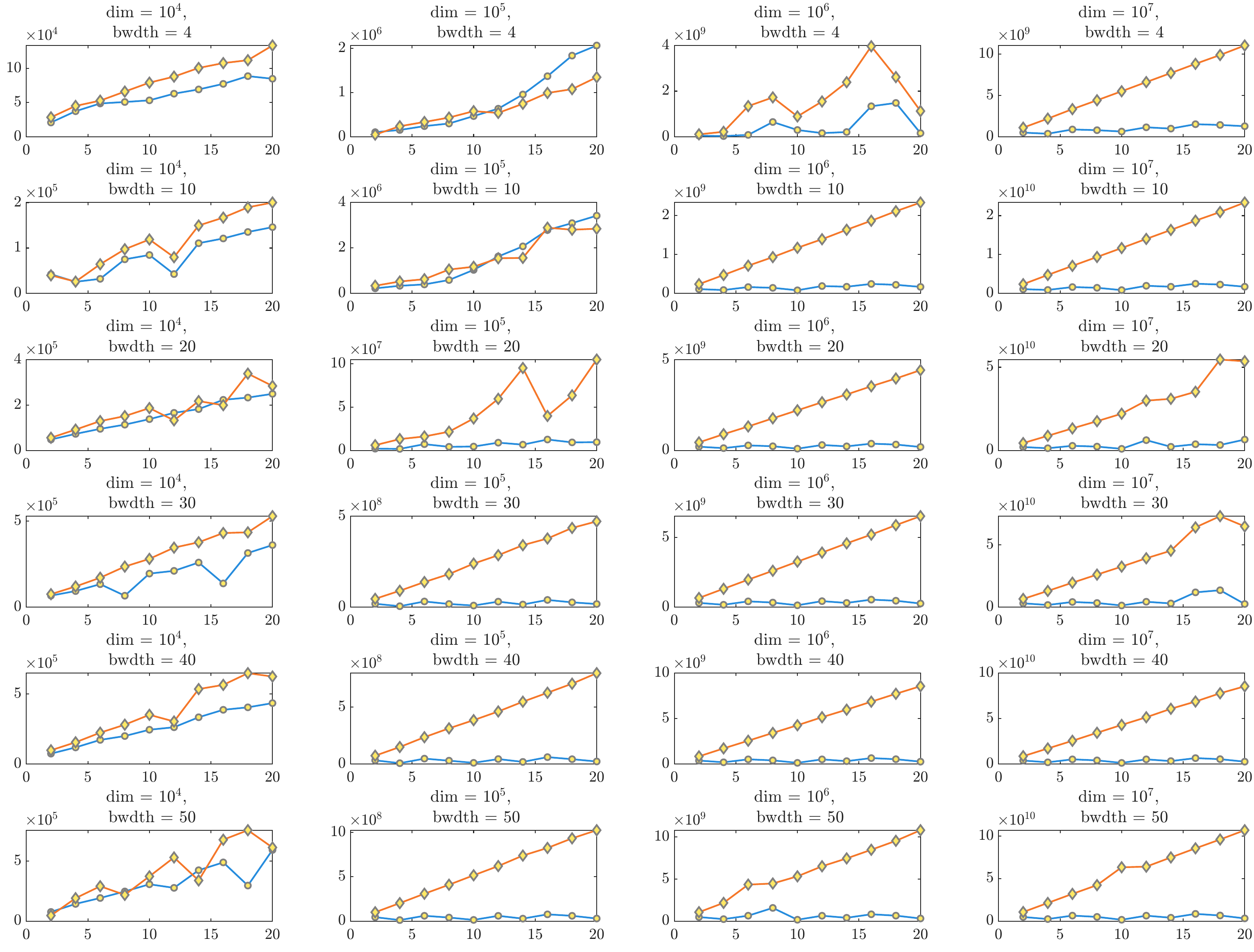}
\caption{\label{figure.cacheMissesSpDiags} For banded matrices of different
dimensions, a comparison of cache misses for the matrix applied to different
sizes of vector blocks both all-at-once (\textbf{blue lines with yellow circles}) and
one-at-a-time (\textbf{red lines with yellow diamonds}).  The horizontal axes show block size and
the vertical axes show number of cache misses per 100 executions of the
experiment. }
\end{figure}

\begin{figure}[htb]
\includegraphics[scale=0.25]{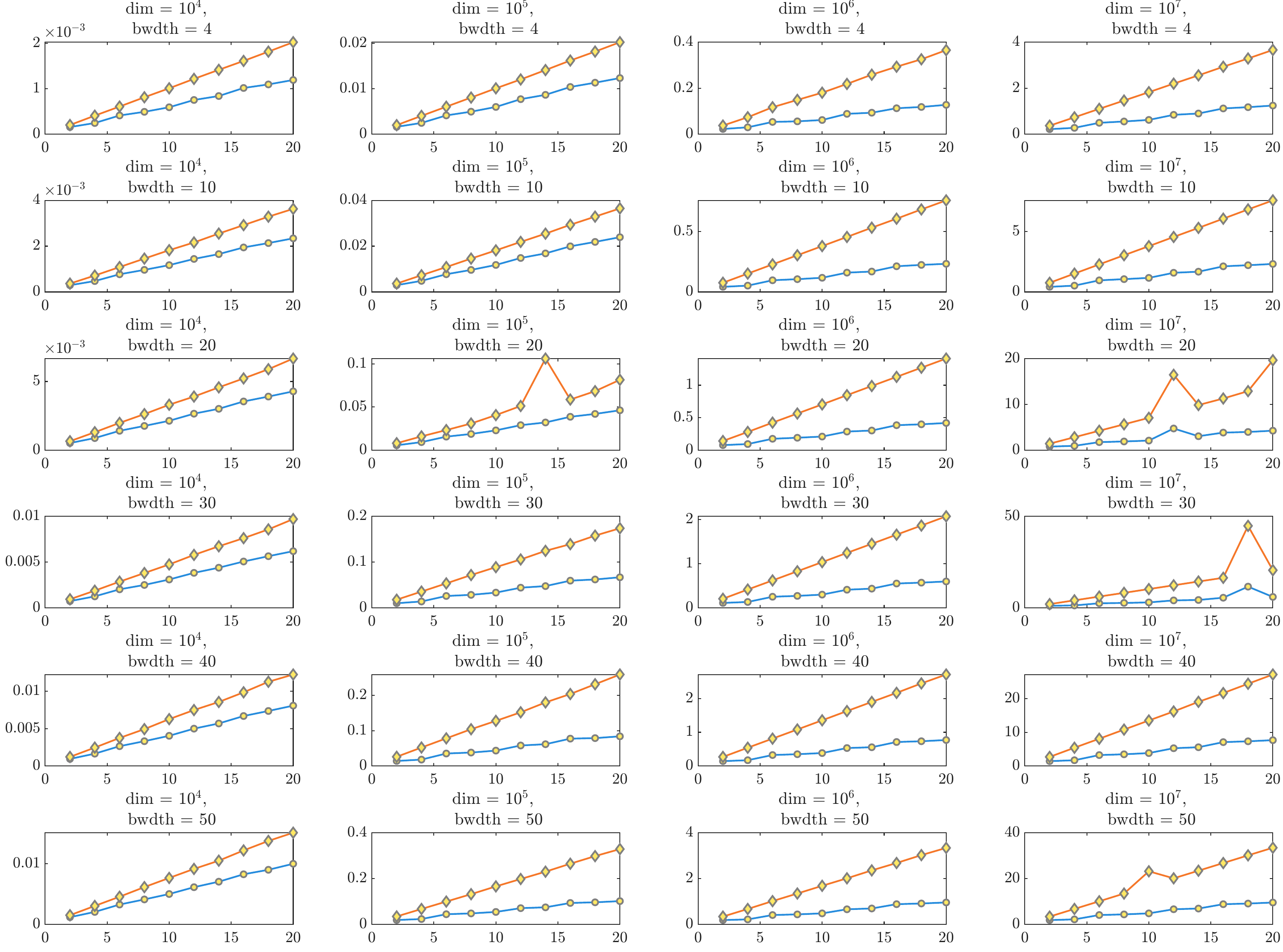}
\caption{\label{figure.avgTimingsSpDiags} For narrow bandwidth matrices of
different dimensions and bandwidths, a comparison average timings (in seconds)
for the projected matrix applied to different sizes of vector blocks both
all-at-once {\rm\textbf{(blue lines with yellow circles)}} and one-at-a-time {\rm\textbf{(gray
lines)}}.  The horizontal axes show block size and the vertical axes show the
average time in seconds required for each experiment. }
\end{figure}

\begin{figure}[htb]
\includegraphics[scale=0.28]{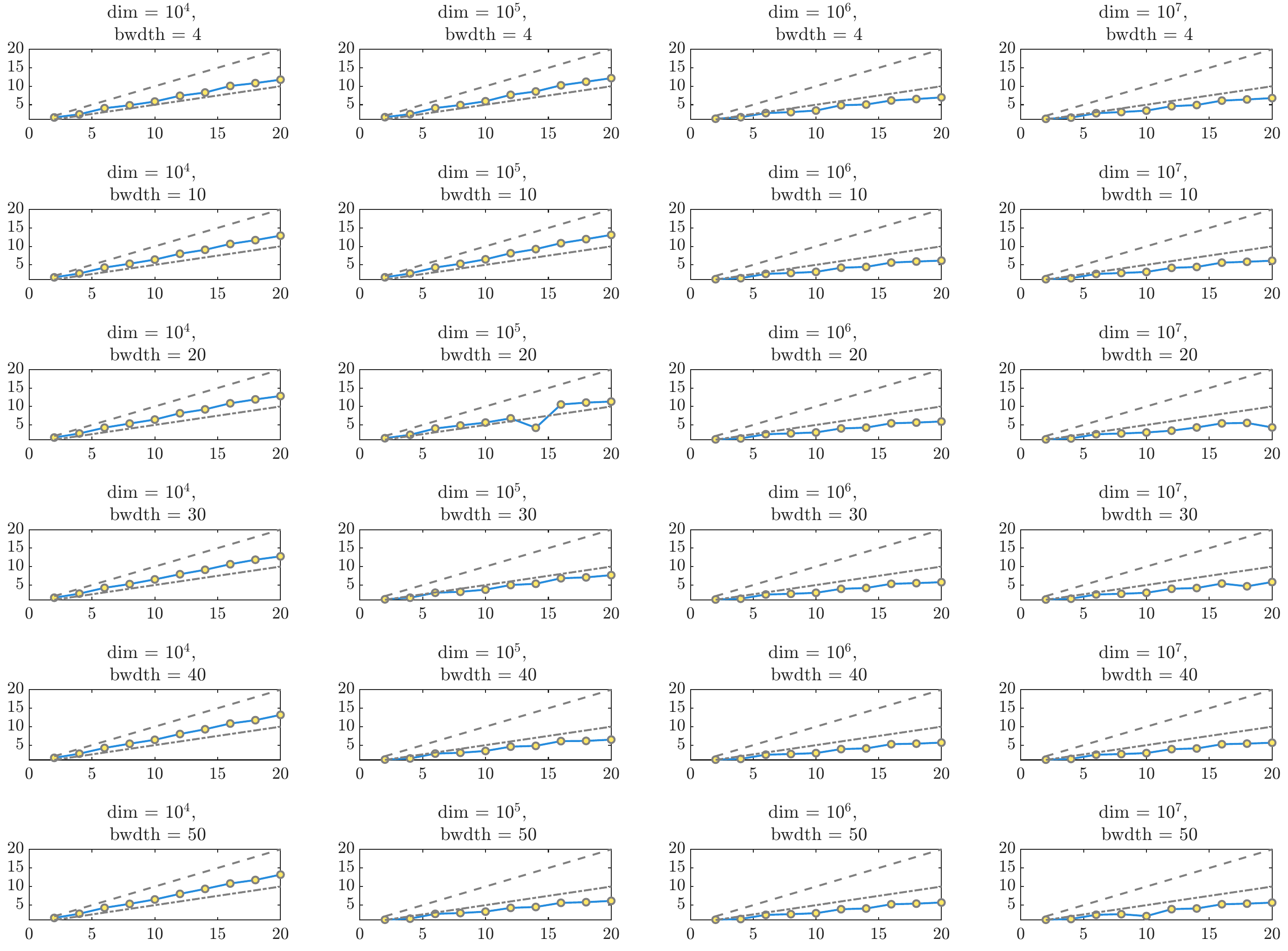}
\caption{\label{figure.blockSingleRatioSpDiags} For narrow bandwidth matrices
of different dimensions and bandwidths, the ratio of the average time taken to
apply the matrix to different sizes of vector blocks versus to a single vector
{\rm\textbf{ (blue lines with yellow circles)}}. The horizontal axes show block size and the
vertical axes show the ratio for each block.  For reference, we include
references lines for ratio $p$ {\rm\textbf{ (gray dashed line)}} and
$\frac{1}{2}p$ {\rm\textbf{ (gray dot-dashed line)}}. }
\end{figure}

In our first set of experiments, we take measurements for fifteen sample matrices arising in a variety of applications, downloaded
from \cite{DH.2011}.  We begin by taking measurements for just the application of the matrix to various sizes of block vectors. In
\Cref{table.uniFloridaMatInfo}, names and relevant characteristics of the matrices are presented.  In
\Cref{figure.cacheMissesFloridaComp}, comparisons of cache misses for single- and block-matvecs are shown for block sizes between
$2$ and $20$.  In \Cref{figure.avgTimingsFloridaComp}, average timings are shown for the same experiments. In
\Cref{figure.blockSingleRatioFloridaComp}, we compare the ratio of the time take to multiply the matrix times a block of vectors
to the time taken to multiply times just one vector. This demonstrates that it is often the case that multiplying times $L$
vectors is not $L$-times as expensive as multiplying times a single vector. In these experiments, we see the greatest
computational benefit for larger matrices, which is when a matvec becomes an I/O-bound operation, i.e., the rate at which data is
used is faster than the rate at which it is retrieved. If one compares the matrix sizes from \Cref{table.uniFloridaMatInfo}
against the data in Figures \ref{figure.cacheMissesFloridaComp}--\ref{figure.blockSingleRatioFloridaComp}, one sees that the
greatest performance difference in the two experiments is for the matrices with the largest number of nonzeros (\texttt{
ML\_Geer}, \texttt{ CoupCons3D}, etc.). For small matrices (e.g., \texttt{epb3}) one observes hardly any difference. This
confirms that using block operations would likely only provides benefit for large matrices. We note that we can transfer these
results to the parallel setting, where we instead consider the situation that the part of the matrix stored on a specific node is
large with respect to the L3 cache size. Note that \Cref{figure.cacheMissesSpDiags} illustrates this relationship;  the figures
further down and to the right show increasingly larger differences in cache misses between the two experiments.

We then took the same measurements but for the projected operator $\prn{\bI - \bPhi}\bA$. We show below experiments for the case
that $\dim\cC = 10$. We also performed the same experiments for $\dim\cC = 50$ and $\dim\cC = 100$, but the results were not
substantially different from those presented.  We see in all three experiments that the cache efficiencies observed for applying the
matrix to blocks of vectors is diminished when a projector is composed with the operator.  In Figures
\ref{figure.cacheMissesProjectedFloridaComp} and \ref{figure.avgTimingsProjectedFloridaComp}, we see that for many matrices, there
is similar performance, in terms of timings and numbers of cache misses, whether the matrix is applied to the full block or to
each vector in the block individually.  In some cases, one-at-a-time application is actually superior. We also again show the
ratio between time taken to multiply the projected matrix times a block of vectors versus multiplying times just one vector.  In
this experiment, we investigate the difference between using modified Gram-Schmidt or multiple passes of classical Gram-Schmidt
(DGKS) \cite{DGKS.1976}.

This experience is important when considering the performance of a block recycled \gmres method as compared to standard block
\gmres.  The cache-efficiency benefits of block methods does not always extend to the orthogonalization routines tested in this
paper.  Thus it may be more appropriate to compare the performance of block \gmres and recycled block \gmres for total search
space dimension being approximately the same, so that the number of orthogonalizations is equivalent for both methods.  There is
much work exploring other methods for efficient, stable orthogonalization in the HPC setting; see e.g.,
\cite{Thomas.Baker.Gaudreault.SISC.25,AudibertGirardonHaddarJolivet:2023:1,CarsonLundRozloznik:2021:1,CarsonMa:2024:1,CarsonLundMaOktay:2024:1,CarsonLundMaOktay:2025:1,JolivetTournier:2016:1}.
Additionally, using all manner of sketched bases to reduce the computational issues concerning orthogonalization is an active area
of exploration \cite{BurkeGuettelSoodhalter:2023:1}, and some authors also explore schemes that do not follow the \gcro approach
and instead choose an unprojected Krylov subspace (i.e., generated by $\bA$ rather than $\prn{\bI - \bPhi}\bA$); see
\cite{BurkeGuettelSoodhalter:2023:1,BurkeFrommerRamirezHidalgoSoodhalter:2022:1}.

In our second set of experiments, we repeat the same tests but for some large, sparse, banded matrices.  These were constructed in
\matlab using \texttt{ rand()} and \texttt{ spdiags() } and saved to disk and then loaded by our compiled code for tests in
\trilinos. The purpose of running tests on such matrices is to
give a clear picture of the cache performance for sparse block $p$ matvecs for matrices whose structure is easily understood but
similar to what often arises in the discretization of differential operators.  These experiments supplement the first set, in
which we use matrices arising from real applications but whose structure is not as simple. We constructed matrices with dimensions
of $10^{4}$, $10^{5}$, $10^{6}$, and $10^{7}$.  Matrices with bandwidths of $4$, $10$, $20$, $30$, $40$, and $50$ were
constructed.  Block sizes tested were the even integers in the interval $\brak{2,20}$.  As the previous experiments demonstrated
that there is a great loss of cache efficiency when applying the projected operator, we restrict these experiments to the case of
the unprojected operator.  In \Cref{figure.cacheMissesSpDiags}, we compare the number of cache misses encountered when applying
the matrices to a block of vector at once versus to the same block one column at-a-time.  In \Cref{figure.avgTimingsSpDiags} we
compare average timings for the same two cases.  In \Cref{figure.blockSingleRatioSpDiags}, we calculate the ratio between the
average time taking to apply each operator to a block of vectors versus applying it to a single vector.  This to a large extent
indicates how much more expensive an iteration (dominated by the cost of the sparse matrix application) will be for a block method
versus the non-block version of that method.  We see again for these artificially constructed sparse matrices that we benefit from
data movement efficiency for block methods.  Furthermore, because the structures of these matrices is precisely known, it is
easier to compare results for different matrices in this group.

\section{Discussion and conclusions}\label{section.concl}
We chose to restrict the experiments in \Cref{section.data-movement} to the matrix application, as it is the main computational
kernel of a Krylov subspace-based iteration. We focus on the performance of the application of $\bA$ and of $(\bI-\bPhi)$ as they
are two of the most dominant costs in Algorithm \ref{alg.block-GCRODR}.  Many of the other operations  are dense matrix-matrix
operations with already confirmed cache efficiency characteristics.  In particular, we use the Householder transformation block
storage and application strategy of Gutknecht and Schmelzer \cite{Gutknecht2008} which means application of the previous
Householder transformations is also a matrix-matrix operation.   Thus the experiments presented in \Cref{section.data-movement}
are analogous to per iteration costs of our method, and for the unprojected operator any block method; see, e.g., \cite{BDJ.2006}.
We have seen that there is a drop in cache efficiency when applying the projected operator.  This is perhaps not surprising.  When
the operator $\bA$ is being applied, perfect cache efficiency arises from being able to store the entirety of $\bA$ (a sparse
matrix) in cache. For the projected operator, we would in addition need to be able to fit $\bC$ (a dense block of vectors) into
cache, as well. This becomes less likely as the number of non-zeros of $\bC$ approaches the number of non-zeros of $\bA$.   The
more vectors we recycle, the less likely that is to happen.  We refer the reader to 
\cite{AudibertGirardonHaddarJolivet:2023:1,JolivetTournier:2016:1} wherein the authors present implementations remedying issues
such as this.

A C++ implementation of Algorithm \ref{alg.block-GCRODR} is available as a part of the \belos package of \trilinos
\cite{PS.GCRODR-Code.2011}, and a \matlab implementation is available at \cite{GCRODR-matlab-code.2016}.  Furthermore, as we have
only shown a small subset of the data coming from our performance results, we also provide tables and \textrm{csv} files
containing the full, raw performance results as a supplement to this paper \cite{RawData_BlockGCRO}.

\section*{Acknowledgments}
The authors would like to thank the referee for the very helpful comments which helped
us improve the presentation. They also thank Florian Tischler of the Johann Radon Institute for Computational and Applied Mathematics for
technical assistance and for allowing the second author to have access to sensitive processor-level functions to take necessary
measurements in \Cref{section.data-movement}. The initial work for this project was undertaken while the second author completed a
graduate research internship at Sandia Laboratories in Albuquerque as a guest of the first author.  The first and second authors
would like to thank Mark Hoemmen for engaging in many helpful conversations and the first author also would like to thank Mike
Heroux for the same.

\textbf{Notice}: This manuscript has been authored by UT-Battelle, LLC, under contract DE-AC05-00OR22725 with the US Department of
Energy (DOE). The US government retains and the publisher, by accepting the article for publication, acknowledges that the US
government retains a nonexclusive, paid-up, irrevocable, worldwide license to publish or reproduce the published form of this
manuscript, or allow others to do so, for US government purposes. DOE will provide public access to these results of federally
sponsored research in accordance with the DOE Public Access Plan (\url{https://www.energy.gov/doe-public-access-plan}).